\documentclass[a4paper]{amsart}
\usepackage{amssymb}

\begin{document}

\def\uL(#1){\L^{\le #1}}

\def\CC{\mathbb C}
\def\NN{\mathbb N}
\def\RR{\mathbb R}
\def\TT{\mathbb T}
\def\ZZ{\mathbb Z}

\def\Bb{{\mathcal B}}
\def\Ee{{\mathcal E}}
\def\Ff{{\mathcal F}}
\def\Hh{{\mathcal H}}
\def\Kk{{\mathcal K}}
\def\Mm{{\mathcal M}}
\def\Oo{{\mathcal O}}
\def\Ss{{\mathcal S}}
\def\Tt{{\mathcal T}}

\def\Ad{\operatorname{Ad}}
\def\Aut{\operatorname{Aut}}
\def\clsp{\overline{\lsp}}
\def\End{\operatorname{End}}
\def\id{\operatorname{id}}
\def\Isom{\operatorname{Isom}}
\def\Ker{\operatorname{Ker}}
\def\lsp{\operatorname{span}}
\def\fin{\operatorname{fin}}
\def\Hom{\operatorname{Hom}}
\def\Obj{\operatorname{Obj}}
\def\dom{\operatorname{dom}}
\def\cod{\operatorname{cod}}
\def\FE{\operatorname{FE}}
\def\Ext{\operatorname{Ext}}
\def\MCE{\operatorname{MCE}}
\def\Mor{\operatorname{Mor}}

\def\Lmin(#1,#2){{%
\L^{\min}(#1,#2)
}}%

\def\L{\Lambda}

\def\SHS{\operatorname{SH\raise1pt\hbox{\hskip-.1em$\times$\hskip-.1em}S}}\def\Ipairs{\SHS(\Lambda)}

\theoremstyle{plain}
\newtheorem{theorem}{Theorem}[section]
\newtheorem*{theorem*}{Theorem}
\newtheorem*{prop*}{Proposition}
\newtheorem{cor}[theorem]{Corollary}
\newtheorem{lemma}[theorem]{Lemma}
\newtheorem{prop}[theorem]{Proposition}
\newtheorem{conj}[theorem]{Conjecture}
\newtheorem{claim}[theorem]{Claim}
\theoremstyle{remark}
\newtheorem{rmk}[theorem]{Remark}
\newtheorem{rmks}[theorem]{Remarks}
\newtheorem*{aside}{Aside}
\newtheorem*{note}{Note}
\newtheorem{comment}[theorem]{Comment}
\newtheorem{example}[theorem]{Example}
\newtheorem{examples}[theorem]{Examples}
\theoremstyle{definition}
\newtheorem{dfn}[theorem]{Definition}
\newtheorem{dfns}[theorem]{Definitions}
\newtheorem{notation}[theorem]{Notation}

\numberwithin{equation}{section}

\title[Ideals in $k$-graph algebras]
{Gauge-invariant ideals in the $C^*$-algebras of finitely aligned higher-rank graphs}
\author{Aidan Sims}
\address{School of Mathematical and Physical Sciences \\ University of Newcastle \\ Callaghan \\ NSW 2308 \\ AUSTRALIA}
\email{Aidan.Sims@newcastle.edu.au}
\keywords{Graphs as categories, graph algebra, $C^*$-algebra}
\date{June 29, 2004}
\subjclass{Primary 46L05}
\thanks{This research is part of the author's Ph.D. thesis, supervised by Professor Iain Raeburn, and was supported by an Australian Postgraduate Award and by the Australian Research Council.}

\begin{abstract}
We produce a complete descrption of the lattice of gauge-invariant ideals in $C^*(\L)$ for a finitely aligned $k$-graph $\L$. We provide a condition on $\L$ under which every ideal is gauge-invariant. We give conditions on $\L$ under which $C^*(\L)$ satisfies the hypotheses of the Kirchberg-Phillips classification theorem.
\end{abstract}

\maketitle

\section{Introduction}\label{sec:intro}
Among the main reasons for the sustained interest in the $C^*$-algebras of directed graphs and their analogues in recent years are the elementary graph-theoretic conditions under which the associated $C^*$-algebra is simple and purely infinite, and the relationship between the gauge-invariant ideals in a graph $C^*$-algebra and the connectivity properties of the underlying graph.

A complete description of the lattice of gauge-invariant ideals of the $C^*$-algebra $C^*(E)$ of a directed graph $E$ was given in \cite{BHRS}, and conditions on $E$ were described under which $C^*(E)$ is simple and purely infinite. Building upon these results, Hong and Szyma\'nski achieved a description of the primitive ideal space of $C^*(E)$ in \cite{HS}. The results of \cite{BHRS} were obtained by a process which builds from a graph $E$ and a gauge-invariant ideal $I$ in $C^*(E)$,  a new graph $F = F(E,I)$ in such a way that the graph $C^*$-algebra $C^*(F)$ is canonically isomorphic to the quotient algebra $C^*(E)/I$. However, recent work of Muhly and Tomforde shows that the quotient algebra $C^*(E)$ can also be regarded as a \emph{relative graph algebra} associated to a subgraph of $E$.

In this note, we turn our attention to the classification of the gauge-invariant ideals in the $C^*$-algebra of a finitely aligned higher-rank graph $\L$, and to the formulation of conditions under which these algebras are simple and purely infinite. Because of the combinatorial peculiarities of higher-rank graphs, constructive methods such as those employed in \cite{BHRS} are not readily available to us in this setting. However, the author has studied a class of \emph{relative Cuntz-Krieger algebras} associated to a higher-rank graph $\L$ in \cite{Si1}, and we use these results to analyse the gauge-invariant ideal structure of $C^*(\L)$. We use the results of \cite{Si1} to give conditions on $\L$ under which $C^*(\L)$ is simple and purely infinite; we also show that relative graph algebras $C^*(\L;\Ee)$, and in particular graph algebras $C^*(\L)$ always belong to the bootstrap class $\mathcal{N}$ of \cite{RSc}, and hence are nuclear and satisfy the UCT.

We begin in Section~\ref{sec:k-graphs} by defining higher-rank graphs, and supplying the definitions and notation we will need for the remainder of the paper. In Section~\ref{sec:hereditary subsets}, we introduce the appropriate analogue in the setting of higher-rank graphs of a \emph{saturated hereditary} set of the vertices of $\L$, and show that such sets $H$ give rise to gauge-invariant ideals $I_H$ in $C^*(\L)$. In Section~\ref{sec:ideals and quotients}, we use the gauge-invariant uniqueness theorem of \cite{Si1} to show that the quotient $C^*(\L)/I_H$ of $C^*(\L)$ by the gauge-invariant ideal associated to a saturated hereditary set $H$ is canonically isomorphic to a relative Cuntz-Krieger algebra $C^*(\L\setminus\L H; \Ee_H)$ associated to a subgraph of $\L$. Using this result, we show in Section~\ref{sec:ideal listing} that the gauge-invariant ideals of $C^*(\L)$ are in bijective correspondence with pairs $(H, B)$ where $H$ is saturated and hereditary, and $B \cup \Ee_H$ is \emph{satiated} as in \cite[Definition~4.1]{Si1}. In Section~\ref{sec:lattice}, we describe the lattice order $\preceq$ on pairs $(H,B)$ which corresponds to the lattice order $\subset$ on gauge-invariant ideals of $C^*(\L)$. In Section~\ref{sec:CK ideals}, we prove that for a certain class of higher-rank graphs $\L$, all the ideals of $C^*(\L)$ are gauge-invariant; however, whilst this result does generalise similar results of \cite{BPRS, RSY1}, the condition~(D) which we need to impose on $\L$ to guarantee that all ideals are gauge-invariant is, in most instances, more or less uncheckable --- the situation is not particulary satisfactory in this regard. In Section~\ref{ch:classifiable} we show that $C^*(\L)$ always falls into the bootstrap class $\mathcal{N}$ of \cite{RSc}, and provide graph-theoretic conditions under which $C^*(\L)$ is simple and purely infinite.

Warning: for consistency with \cite{KP}, the author has continued to use terminology such as ``hereditary'' and ``cofinal'' in this paper. Readers familiar with graph algebras should be wary as to the meaning of these terms because of the change of edge-direction conventions involved in going from directed graphs to $k$-graphs.

{\bf Acknowledgements.} This article is based on part of the author's PhD dissertation, which was written at the University of Newcastle, Australia, under the supervision of Iain Raeburn. The author would like to thank Iain for his insight and guidance. The author would also like to thank D. Gwion Evans for directing his attention to the results of \cite{PQR} on AF algebras associated to non-row-finite skew-product $k$-graphs.

\section{Higher-rank graphs and their representations} \label{sec:k-graphs}
The definitions in this section are taken more or less wholesale from \cite{Si1}.

We regard $\NN^k$ as an additive semigroup with identity 0. For $m,n \in \NN^k$, we write $m \vee n$ for their coordinate-wise maximum and $m \wedge n$ for their coordinate-wise minimum. We write $n_i$ for the $i^{\rm th}$ coordinate of $n \in \NN^k$, and $e_i$ for the $i^{\rm th}$ generator of $\NN^k$; so $n = \sum^k_{i=1} n_i \cdot e_i$.

\begin{dfn} \label{dfn:k-graph} 
Let $k \in \NN \setminus \{0\}$.  A $k$-graph is a pair $(\L,d)$ where $\L$ is a countable category and $d$ is a functor from $\L$ to $\NN^k$ which satisfies the \emph{factorisation property}:  {\sl For all $\lambda \in \Mor(\L)$ and all $m,n \in \NN^k$ such that $d(\lambda) = m+n$, there exist unique morphisms $\mu$ and $\nu$ in $\Mor(\L)$ such that $d(\mu) = m$, $d(\nu) = n$ and $\lambda = \mu\nu$.}
\end{dfn}

Since we are regarding $k$-graphs as generalised graphs, we refer to elements of $\Mor(\L)$ as \emph{paths} and we write $r$ and $s$ for the codomain and domain maps.

The factorisation property implies that $d(\lambda) = 0$ if and only if $\lambda = \id_v$ for some $v \in \Obj(\L)$. Hence we identify $\Obj(\L)$ with $\{\lambda \in \Mor(\L) : d(\lambda) = 0\}$, and write $\lambda \in \L$ in place of $\lambda \in \Mor(\L)$.

Given $\lambda \in \L$ and $E \subset\L$, we define $\lambda E :=\{\lambda\mu : \mu \in E, r(\mu) = s(\lambda)\}$ and $E\lambda := \{\mu\lambda : \mu \in E, s(\mu) = r(\lambda)\}$.  In particular if $d(v) = 0$, then $vE = \{\lambda \in E : r(\lambda) = v\}$.  In analogy with the path-space notation for $1$-graphs, we denote by $\L^n$ the collecton $\{\lambda \in \L : d(\lambda) = n\}$ of paths of degree $n$ in $\L$.

The factorisation property ensures that if $l \le m \le n \in \NN^k$ and if $d(\lambda) = n$, then there exist unique elements, denoted $\lambda(0, l)$, $\lambda(l, m)$ and $\lambda(m,n)$, of $\L$ such that $d(\lambda(0,l)) = l$, $d(\lambda(l,m)) = m-l$, and $d(\lambda(m,n)) = n-m$ and such that $\lambda = \lambda(0,l)\lambda(l,m)\lambda(m,n)$.

\begin{dfn} \label{dfn:common extensions} 
Let $(\L,d)$ be a $k$-graph. For $\mu,\nu \in \L$ we denote the collection $\{\lambda \in \L : d(\lambda) = d(\mu) \vee d(\nu), \lambda(0, d(\mu)) = \mu, \lambda(0, d(\nu)) = \nu\}$ of \emph{minimal common extensions} of $\mu$ and $\nu$ by $\MCE(\mu,\nu)$. We write $\Lmin(\mu,\nu)$ for the collection
\[
\Lmin(\mu,\nu) := \{(\alpha,\beta) \in \L \times \L : \mu\alpha =  \nu\beta \in \MCE(\mu,\nu)\}.
\]
If $E \subset \L$ and $\mu \in \L$, then we write $\Ext_\L(\mu;E)$ for the set
\[\begin{split}
\Ext_\L(\mu;E) := \{\beta \in s(\mu)\L : \text{there }&\text{exists } \nu \in E \text{ such that } \mu\beta \in \MCE(\mu,\nu)\};
\end{split}\]
when the ambient $k$-graph $\L$ is clear from context, we write $\Ext(\mu;E)$ in place of $\Ext_\L(\mu; E)$. We say that $\L$ is finitely aligned if $|\MCE(\mu,\nu)| < \infty$ for all $\mu,\nu \in \L$. 

Let $v \in \L^0$ and $E \subset v\L$. We say $E$ is \emph{exhaustive} if $\Ext(\lambda;E) \not= \emptyset$ for all $\lambda \in v\L$.
\end{dfn}

\begin{notation} 
Let $(\L,d)$ be a finitely aligned $k$-graph. Define
\[
\FE(\L) := \textstyle\bigcup_{v \in \L^0} \{E \subset v\L\setminus\{v\} :\text{ $E$ is finite and exhaustive}\}.
\]
For $E \in \FE(\L)$ we write $r(E)$ for the vertex $v \in \L^0$ such that $E \subset v\L$.
\end{notation}

Notice that whilst any finite subset of $v\L$ which contains $v$ is automatically finite exhaustive, we don not include such sets in $\FE(\L)$. Note also that since $v\L$ is never empty (in particular, it always contains $v$), finite exhausitve sets, and in particular elements of $\FE(\L)$, are always nonempty.

\begin{dfn}\label{dfn:relCK family}
Let $(\L,d)$ be a finitely aligned $k$-graph, and let $\Ee$ be a subset of $\FE(\L)$. A \emph{relative Cuntz-Krieger $(\L;\Ee)$-family} is a collection $\{t_\lambda : \lambda \in \L\}$ of partial isometries in a $C^*$-algebra satisfying
\begin{itemize}
\item[(TCK1)] $\{t_v : v \in \L^0\}$ is a collection of mutually orthogonal projections;
\item[(TCK2)] $t_\lambda t_\mu = \delta_{s(\lambda), r(\mu)} t_{\lambda\mu}$ for all $\lambda, \mu \in \L$;
\item[(TCK3)] $t^*_\lambda t_\mu = \sum_{(\alpha,\beta) \in \Lmin(\lambda,\mu)} t_\alpha t^*_\beta$ for all $\lambda,\mu \in \L$; and
\item[(CK)] $\prod_{\lambda \in E} (t_{r(E)} - t_\lambda t^*_\lambda) = 0$ for all $E \in \Ee$.
\end{itemize}
When $\Ee = \FE(\L)$, we call $\{t_\lambda : \lambda \in \L\}$ a \emph{Cuntz-Krieger $\L$-family}.
\end{dfn}

For each pair $(\L,\Ee)$ there exists a universal  $C^*$-algebra $C^*(\L;\Ee)$, generated by a universal relative Cuntz-Krieger $(\L;\Ee)$-family $\{s_\Ee(\lambda) : \lambda \in \L\}$ which admits a \emph{gauge-action} $\gamma$ of $\TT^k$ satisfying $\gamma_z(s_\Ee(\lambda)) = z^{d(\lambda)} s_\Ee(\lambda)$. We write $C^*(\L)$ for $C^*(\L;\FE(\L))$, and call it the \emph{Cuntz-Krieger} algebra, and we denote the universal Cuntz-Krieger family by $\{s_\lambda : \lambda \in \L\}$; this agrees with the definitions given in \cite{RSY2}.

There is also a Toeplitz algebra $\Tt C^*(\L)$ associated to each $k$-graph $\L$. By definition, this is the universal $C^*$-algebra generated by a family $\{s_{\Tt}(\lambda) : \lambda \in \L\}$ which satisfy (TCK1)--(TCK3), and hence is canonically isomorphic to $C^*(\L;\emptyset)$. Indeed, each $C^*(\L;\Ee)$ is a quotient of $\Tt C^*(\L)$:

\begin{lemma}\label{lem:univ-quotients}
Let $(\L,d)$ be a finitely aligned $k$-graph, and let $\Ee \subset \FE(\L)$. Let $J_{\Ee}$ denote the ideal of $\Tt C^*(\L)$ generated by the projections 
\[
\textstyle
\Big\{\prod_{\lambda \in E} \big(s_\Tt(r(E)) - s_\Tt(\lambda) s_\Tt(\lambda)^*\big) : E \in \Ee\Big\}.
\]
Then $C^*(\L;\Ee)$ is canonically isomorphic to $\Tt C^*(\L) / J_\Ee$.
\end{lemma}
\begin{proof}
The universal property of $\Tt C^*(\L)$ gives a homomorphism $\pi : \Tt C^*(\L) \to C^*(\L;\Ee)$ satisfying $\pi(s_\Tt(\lambda)) = s_\Ee(\lambda)$ for all $\lambda$. Since $\{s_\Ee(\lambda) : \lambda \in \L\}$ satisfy~(CK), we have $J_\Ee \subset \ker\pi$ and hence $\pi$ descends to a homomorphism $\tilde\pi : \Tt C^*(\L) / J_\Ee \to C^*(\L;\Ee)$ such that $\tilde\pi(s_\Tt(\lambda) + J_\Ee) = s_\Ee(\lambda)$ for all $\lambda$.

On the other hand, the family $\{s_\Tt(\lambda) + J_\Ee : \lambda \in \L\} \subset \Tt C^*(\L) / J_{\Ee}$ satisfy~(CK) by definition of $J_\Ee$, so the universal property of $C^*(\L;\Ee)$ gives a homomorphism $\phi : C^*(\L;\Ee) \to \Tt C^*(\L) / J_{\Ee}$ such that $\phi(s_\Ee(\lambda)) = s_\Tt(\lambda) + J_\Ee$ for all $\lambda$. We have that $\tilde\pi$ and $\phi$ are mutually inverse, and the result follows.
\end{proof}

\section{Hereditary subsets and associated ideals} \label{sec:hereditary subsets}

\begin{dfn}\label{dfn:sat,hered} Let $(\L,d)$ be a finitely aligned $k$-graph.  Define a relation $\le$ on $\L^0$ by $v \le w$ if and only if $v \L w
\not= \emptyset$.
\begin{itemize}
\item[(1)] We say that a subset $H$ of $\L^0$ is \emph{hereditary} if $v \in H$ and $v \le w$ imply $w \in H$.
\item[(2)] We say that $H \subset \L^0$ is \emph{saturated} if, whenever $v \in \L^0$ and there exists a finite exhaustive subset $F \subset v \L$ with $s(F) \subset H$, we also have $v \in H$.
\end{itemize} 
For $H \subset \L^0$ we call the smallest saturated set containing $H$ the \emph{saturation} of $H$.
\end{dfn}

\begin{lemma} \label{lem:saturation} 
Let $(\L,d)$ be a finitely aligned $k$-graph and let $G \subset \L^0$. Let $\Sigma G := \{v \in \L^0 : \text{ there exists a finite exhaustive set } F \subset v\L G\}$. Then
\begin{itemize}
\item[(1)] $\Sigma G$ is equal to the saturation of $G$; and
\item[(2)] if $G$ is hereditary, then $\Sigma G$ is hereditary.
\end{itemize}
\end{lemma}
\begin{proof} 
First note that if $v \in G$ then $\{v\} \subset v\L G$ is finite and exhaustive so that $G \subset \Sigma G$. Note also that $\Sigma G$ is a subset of the saturation of $G$ by definition. To see that $\Sigma G$ is saturated, let $v \in \L^0$ and suppose $F \in v\L(\Sigma G)$ is finite and exhaustive. If $v \in F$, then $v \in \Sigma G $ by definition, so suppose that $v \not\in F$. Let $E := \{\lambda \in F : s(\lambda) \not \in G\}$. By definition of $\Sigma G $, for each $\lambda \in E$, there exists $E_\lambda \in s(\lambda)\FE(\L)$ with $s(E_\lambda) \subset G$. Then \cite[Lemma~5.3]{Si1} shows that $F' := (F \setminus E) \cup (\bigcup_{\lambda \in E} \lambda E_\lambda)$ belongs to $\FE(\L)$. Since $F' \subset v\L G$, it follows that $v \in \Sigma G $ by definition. This establishes~(1).

To prove claim (2), suppose $G$ is hereditary, and suppose $v,w \in \L^0$ satisfy $v \in \Sigma G$ and $v \le w$; say $\lambda \in \L$ with $r(\L) = v$, $s(\L) = w$.  If $v \in G$ then $w \in G$ because $G$ is hereditary, so suppose that $v \in \Sigma G  \setminus G$.  By definition of $\Sigma$ there exists $F \in v \FE(\L)$ such that $s(F) \subset G$.  By \cite[Lemma~2.3]{Si1}, $\Ext(\lambda;F)$ is a finite exhaustive subset of $w \L$. Since $s(F) \subset G$, and since, for $\alpha \in \Ext(\lambda;F)$, we have $s(\alpha) \le s(\mu)$ for some $\mu \in F$, we have $s(\Ext(\lambda;F)) \subset G$. It follows that $w \in \Sigma G $, completing the proof.
\end{proof}

\begin{lemma} \label{lem:H_I sat,hered} 
Let $(\L,d)$ be a finitely aligned $k$-graph, and let $I$ be an ideal of $C^*(\L)$.  Then $H_I := \{v \in \L^0 : s_v \in I\}$ is saturated and hereditary.
\end{lemma}

To prove Lemma~\ref{lem:H_I sat,hered}, we first need to recall some notation from \cite{RS1}.

\begin{notation}\label{ntn:vee E}
Let $(\L,d)$ be a finitely aligned $k$-graph and let $E$ be a finite subset of $\L$. As in \cite{RS1}, we denote by $\vee E$ the smallest subset of $\L$ such that $E \subset \vee E$ and such that if $\lambda,\mu \in \vee E$, then $\MCE(\lambda,\mu) \subset \vee E$. We have that $\vee E$ is finite and that $\lambda \in \vee E$ implies $\lambda = \mu\mu'$ for some $\mu \in E$ by \cite[Lemma~8.4]{RS1}. 
\end{notation}

\begin{proof}[Proof of Lemma~\ref{lem:H_I sat,hered}]
Suppose $v \in H_I$ and $w \in \L^0$ with $v \le w$.  So there exists $\lambda \in v\L w$. Since $s_v \in I$, we have $s_w = s_\lambda^* s_v s_\lambda \in I$, and then $w \in H_I$; consequently $H_I$ is hereditary.  Now suppose that $v \in \L^0$ and there is a finite exhaustive set $F \subset v \L$ with $s(F) \subset H_I$.  By \cite[Lemma~3.1]{RSY2}, we have $s_v \in \lsp\{s_\lambda s^*_\lambda : \lambda \in \vee F\}$.  Since $\lambda \in \vee F$ implies $\lambda = \alpha\alpha'$ for some $\alpha \in F$, and since $H_I$ is hereditary, we have $s(\vee F) \subset H_I$.  Consequently, for $\lambda \in \vee F$, we have $s_\lambda s^*_\lambda = s_\lambda s_{s(\lambda)} s^*_\lambda \in I$, so $s_v \in I$, giving $v \in H_I$.
\end{proof}

\begin{notation}\label{ntn:I_H}
For $H \subset \L^0$, let $I_H$ be the ideal in $C^*(\L)$ generated by $\{s_v : v \in H\}$. Let $H\L$ denote the subcategory $\{\lambda \in \L : r(\lambda) \in H\}$ of $\L$.
\end{notation}

\begin{lemma}\label{lem:Lambda(H) Vs I_H} 
Let $(\L,d)$ be a finitely aligned $k$-graph, and suppose that $H \subset \L^0$ is saturated and hereditary.  Then $(H\L, d|_{H\L})$ is also a finitely aligned $k$-graph, and $C^*(H\L) \cong C^*(\{s_\lambda : r(\lambda) \in H\}) \subset C^*(\L)$. Moreover this subalgebra is a full corner in $I_H$.
\end{lemma}
\begin{proof} 
One checks that $(H\L, d|_{H\L})$ is a $k$-graph just as in \cite[Theorem~5.2]{RSY1}, and it is finitely aligned because $(H\L)^{\min}(\lambda,\mu) \subset \Lmin(\lambda,\mu)$.

The universal property of $C^*(H\L)$ ensures that there exists a homomorphism $\pi : C^*(H\L) \to C^*(\{s_\lambda : r(\lambda) \in H\})$. Write $\gamma_H$ for gauge action on $C^*(H\L)$ and $\gamma|$ for the restriction of the gauge action on $C^*(\L)$ to $C^*(\{s_\lambda : r(\lambda) \in H\})$.  Then $\pi \circ (\gamma_H)_z = (\gamma|)_z \circ \pi$ for all $z \in \TT^k$, and \cite[Theorem~4.2]{RSY2} shows that $\pi$ is injective.

For the final statement, just use the argument of \cite[Theorem~4.1(c)]{BPRS} to see that $C^*(\{s_\lambda : r(\lambda) \in H\})$ is the corner of $I_H$ determined by the projection $P_H := \sum_{v \in H} s_v \in \Mm(I_H)$, and that this projection is full.
\end{proof}

\section{Quotients of $C^*(\L)$ by $I_H$} \label{sec:ideals and quotients} 
We now want to show that the quotients of Cuntz-Krieger algebras by the ideals $I_H$ of section~\ref{sec:hereditary subsets} are relative Cuntz-Krieger algebras associated to $\L \setminus \L H$.

Let $(\L,d)$ be a $k$-graph, and let $H \subset \L^0$ be a saturated hereditary set. Consider the subcategory $\L\setminus\L H = \{\lambda \in \L : s(\lambda) \not \in H\}$.

\begin{lemma}\label{lem:quotient graph} 
Let $(\L,d)$ be a finitely aligned $k$-graph, and let $H \subset \L^0$ be saturated and hereditary.  Then $(\L\setminus\L H, d|_{\L\setminus\L H})$ is also a finitely aligned $k$-graph.
\end{lemma}
\begin{proof} 
We first check the factorisation property for $(\L\setminus\L H, d|_{\L\setminus\L H})$, and then that $(\L\setminus\L H, d|_{\L\setminus\L H})$ is finitely aligned.  For the factorisation property, let $\lambda \in \L\setminus\L H$, and let $m,n \in \NN^k$, $m+n = d(\lambda)$.  By the factorisation property for $\L$, there exist unique $\mu,\nu \in \L$ such that $d(\mu) = m$, $d(\nu) = n$ and $\lambda = \mu\nu$.  Since $s(\nu) = s(\lambda) \not \in H$, we have $\nu \in \L\setminus\L H$.  Since, by definition of $\le$, we have $r(\nu) \le s(\nu)$ it follows that $r(\nu) \not\in H$ because $H$ is hereditary.  But $r(\nu) = s(\mu)$ so it follows that $\mu \in \L\setminus\L H$.  Finite alignedness of the $k$-graph $\L\setminus\L H$ is trivial since $(\L\setminus\L H)^{\min}(\lambda,\mu) \subset \Lmin(\lambda,\mu)$ for all $\lambda,\mu \in \L\setminus\L H$.
\end{proof}

\begin{dfn}\label{dfn:EesubH}
Let $(\L,d)$ be a finitely aligned $k$-graph and let $H$ be a saturated hereditary subset of $\L^0$. Define $\Ee_H := \{E \setminus EH : E \in \FE(\L)\}$.
\end{dfn}

\begin{lemma}\label{lem:EeH f.e.} 
Let $(\L,d)$ be a finitely aligned $k$-graph, and suppose that $H \subset \L^0$ is saturated and hereditary.  Then $\Ee_H \subset \FE(\L\setminus \L H)$.
\end{lemma}
\begin{proof} 
Suppose that $E \in \Ee_H$ and that $\mu \in r(E)(\L \setminus \L H)$.  Suppose for contradiction that $(\L \setminus \L H)^{\min}(\lambda,\mu) = \emptyset$ for all $\lambda \in E$.  Since $E \in \Ee_H$, there exists $F \in \FE(\L)$ such that $F \setminus F H = E$.  We have
\begin{equation}
\label{eqn:Exts}
\Ext_\L(\mu; F) = \Ext_\L(\mu;E) \cup \Ext_\L(\mu; F \setminus E) = \Ext_\L(\mu;E) \cup \Ext_\L(\mu; F H).
\end{equation}
Now $F H \subset \L H$ by definition, and then $\Ext(\mu; F H) \in \L H$ because $H$ is hereditary. Since $(\L \setminus \L H)^{\min}(\lambda,\mu) = \emptyset$ for all $\lambda \in E$, we must have $\Lmin(\lambda,\mu) \subset \L H \times \L H$ for all $\lambda \in E$, and hence we also have $\Ext_\L(\mu;E) \subset \L H$.  Hence \eqref{eqn:Exts} shows that $\Ext_\L(\mu;F) \subset \L H$.  But $F$ is exhaustive in $\L$, so $\Ext(\mu;F)$ is also exhaustive by \cite[Lemma~2.3]{Si1}, and then since $H$ is saturated, it follows that $s(\mu) \in H$, contradicting our choice of $\mu$.
\end{proof}

\begin{theorem}\label{thm:quotient conjecture} Let $(\L,d)$ be a finitely aligned $k$-graph, and let $H \subset \L^0$ be saturated and hereditary.  Then $C^*(\L) / I_H$ is canonically isomorphic to $C^*((\L\setminus\L H) ; \Ee_H)$.
\end{theorem}

To prove Theorem~\ref{thm:quotient conjecture}, we need to collect some additional results. Recall from
\cite[Definition~4.1]{Si1} that a subset $\Ee$ of $\FE(\L)$ is said to be \emph{satiated} if it satisfies
\begin{itemize}
\item[(S1)] if $G \in \Ee$ and $E \in \FE(\L)$ with $G \subset E$, then $E \in \Ee$;
\item[(S2)] if $G \in \Ee$ with $r(G) = v$ and $\mu \in v\L \setminus G\L$, then $\Ext(\mu;G) \in \Ee$;
\item[(S3)] if $G \in \Ee$ and $0 < n_\lambda \le d(\lambda)$ for $\lambda \in G$, then $\{\lambda(0, n_\lambda) : \lambda \in G\} \in \Ee$; and
\item[(S4)] if $G \in \Ee$, $G' \subset G$ and for each $\lambda \in G'$, $G'_\lambda$ is an element of $\Ee$ such that $r(G_\lambda') = s(\lambda)$, then $\textstyle\big((G\setminus G') \cup \big(\bigcup_{\lambda \in G'} \lambda
G'_\lambda\big)\big) \in \Ee$.
\end{itemize}

\begin{lemma}\label{lem:EsubH its own bar} Let $(\L,d)$ be a finitely aligned $k$-graph, and let $H \subset \L^0$ be saturated and hereditary. Then $\Ee_H$ is satiated.
\end{lemma}
\begin{proof} 
For (S1), suppose that $E \in \Ee_H$ and $F \subset \L \setminus \L H$ is finite with $E \subset F$.  By definition of $\Ee_H$, there exists $E' \in \FE(\L)$ such that $E' \setminus E' H = E$.  But then $F' := F \cup E' H \in \FE(\L)$ by \cite[Lemma~5.3]{Si1}. Since $F = F' \setminus F' H$, it follows that $F \in \Ee_H$.

For (S2), suppose that $E \in \Ee_H$, that $\mu \in r(E)(\L \setminus \L H)$ and that $\mu \not\in E\L$.  Since $E \in \Ee_H$, there exists $E' \in \FE(\L)$ such that $E' \setminus E' H = E$.  Since $\mu \in \L \setminus \L H$, we have $\mu \not\in E' H$, and hence $\Ext_\L(\mu;E') \in \FE(\L)$ by \cite[Lemma~2.3]{Si1}. We also have
\[
\begin{split}
\Ext_\L(\mu;E') 
&= \Ext_\L(\mu; E) \cup \Ext_\L(\mu; E' H) \\
&= \Ext_{\L \setminus \L H}(\mu; E) \cup \Ext_\L(\mu;E) H \cup \Ext_\L(\mu;E' H).
\end{split}
\]
Since both $\Ext\L(\mu;E) H$ and $\Ext_\L(\mu; E'H)$ are subsets of $\L H$, it follows that
\[
\Ext_{\L \setminus \L H}(\mu;E) = \Ext_\L(\mu;E') \setminus \Ext_\L(\mu;E') H,
\]
and hence belongs to $\Ee_H$. 

For (S3), suppose that $E \in \Ee_H$, say $E' \in \FE(\L)$ and $E = E' \setminus E' H$.  For each $\lambda \in E$, let $n_\lambda \in \NN^k$ with $0 < n_\lambda \le d(\lambda)$.  For $\mu \in E' H$, let $n_\mu := d(\mu)$.  Since $E'$ is exhaustive in $\L$, we have that $\{\mu(0,n_\mu) : \mu \in E'\}$ is also a finite exhaustive subset of $\L$ by \cite[Lemma~5.3]{Si1}, and since
\[
\{\lambda(0,n_\lambda) : \lambda \in E\} = \{\mu(0, n_\mu) : \mu \in E'\} \setminus \{\mu(0, n_\mu) : \mu \in E' H\},
\]
it follows that $\{\lambda(0,n_\lambda) : \lambda \in E\} \in \Ee_H$.

Finally, for (S4), suppose that $E \in \Ee_H$, say $E' \in \FE(\L)$ and $E = E' \setminus E'H$.  Let $F \subset E$, and for each $\lambda \in F$, suppose that $F_\lambda \in \Ee_H$ with $r(F_\lambda) = s(\lambda)$.  We must show that $G := (E \setminus F) \cup \big(\bigcup_{\lambda \in F} \lambda F_\lambda\big) \in \Ee_H$.  Since each $F_\lambda \in \Ee_H$, for each $\lambda \in F$, there exists a set $F'_\lambda \in \FE(\L)$ with $F_\lambda = F'_\lambda \setminus F'_\lambda H$.  Let $G' := (E' \setminus F) \cup \big(\bigcup_{\lambda \in F} \lambda F'_\lambda\big)$.  We will show that $G = G' \setminus G'H$, and that $G'$ is finite and exhaustive in $\L$; it follows from the definition of $\Ee_H$ that $G \in \Ee_H$, proving the result.

We have $G' \in \FE(\L)$ by \cite[Lemma~5.3]{Si1}, so it remains only to show that $G = G' \setminus G'H$.  But since $H$ is hereditary, we have
\[\begin{split}
G'H &=\textstyle \Big((E' \setminus F) \cup \big(\bigcup_{\lambda \in F} \lambda F'_\lambda\big)\Big) H \\
&=\textstyle (E' \setminus F) H \cup \big(\bigcup_{\lambda \in F} \lambda(F'_\lambda H)\big) 
= E'H \cup \big(\bigcup_{\lambda \in F} \lambda F'_\lambda\big) H
\end{split}\]
because $F \subset E \subset \L\setminus\L H$. Consequently
\[\textstyle
G' \setminus G'H = \big((E' \setminus F) \cup \big(\bigcup_{\lambda \in F} \lambda F'_\lambda\big)\big) 
\setminus 
\big(E'H \cup \big(\bigcup_{\lambda \in F} \lambda F'_\lambda H\big)\big) = G
\]
as required.
\end{proof}

\begin{lemma}\label{lem:CK fams descend} 
Let $(\L,d)$ be a finitely aligned $k$-graph, and let $H \subset \L^0$ be saturated and hereditary.  Let $\{t_\lambda :\lambda \in \L\}$ be a Cuntz-Krieger $\L$-family, and let $I^t_H$ be the ideal in $C^*(\{t_\lambda : \lambda \in \L\})$ generated by $\{t_v : v \in H\}$.  Then $\{t_\lambda + I^t_H : \lambda \in \L\setminus\L H\}$ is a relative Cuntz-Krieger $(\L \setminus\L H; \Ee_H)$-family in $C^*(\{t_\lambda : \lambda \in \L\}) / I^t_H$.
\end{lemma}
\begin{proof} 
Relations (TCK1) and (TCK2) hold automatically since they also hold for the Cuntz-Krieger $\L$-family $\{t_\lambda : \lambda \in \L\}$.  For (TCK3), let $\lambda,\mu \in \L \setminus\L H$ and notice that since $\{t_\lambda : \lambda \in \L\}$ is a Cuntz-Krieger $\L$-family, we have
\[
(t^*_\lambda + I^t_H)(t_\mu + I^t_H) = \sum_{(\alpha,\beta) \in \Lmin(\lambda,\mu)} t_\alpha t^*_\beta + I^t_H.
\]
To show that this is equal to $\sum_{(\alpha,\beta) \in (\L\setminus\L H)^{\min}(\lambda,\mu)} t_\alpha t^*_\beta + I^t_H$, we need to show that
\[
(\alpha,\beta) \in \Lmin(\lambda,\mu) \setminus (\L \setminus \L H)^{\min}(\lambda, \mu) \text{ implies } t_\alpha t^*_\beta \in I^t_H.
\]
So fix $(\alpha,\beta) \in \Lmin(\lambda,\mu) \setminus (\L \setminus \L H)^{\min}(\lambda, \mu)$.  Then $s(\alpha) = s(\beta) \in H$, and hence $s_\alpha s^*_\beta = s_\alpha s_{s(\alpha)} s^*_\beta \in I^t_H$.

It remains to check~(CK).  Let $E \in \Ee_H$, say $E' \in \FE(\L)$ and $E = E' \setminus E' H$, and let $v := r(E)$.  We must show that $\prod_{\lambda \in E} (t_v - t_\lambda t^*_\lambda)$ belongs to $I^t_H$.  We know that $\prod_{\lambda \in E'} (t_v - t_\lambda t^*_\lambda) = 0$, and it follows that
\begin{equation}\label{eq:orth to FH}
\textstyle
\prod_{\lambda \in E}(t_v - t_\lambda t^*_\lambda) \Big(\prod_{\mu \in E' H} (t_v - t_\mu t^*_\mu)\Big) = 0.
\end{equation}
Since $H$ is hereditary, Notation~\ref{ntn:vee E} gives $\vee (E'H) \subset \L H$, and $\prod_{\mu \in \vee (E' H)} (t_v - t_\mu t^*_\mu) \le \prod_{\mu \in E' H} (t_v - t_\mu t^*_\mu)$.  Furthermore by \cite[Proposition~3.5]{RSY2} we have
\[
\textstyle
t_v = \prod_{\mu \in \vee (E' H)} (t_v - t_\mu t^*_\mu) + \sum_{\mu \in \vee (E' H)} Q(t)^{\vee (E' H)}_\mu
\]
where $Q(t)^{\vee (E' H)}_\mu := \prod_{\mu\mu' \in \vee(E' H) \setminus \{\mu\}} (t_\mu t^*_\mu - t_{\mu\mu'} t^*_{\mu\mu'})$.

Hence we can calculate
\begin{align*}
\prod_{\lambda \in E}(t_v - t_\lambda t^*_\lambda) 
&= \Big(\prod_{\lambda \in E}(t_v - t_\lambda t^*_\lambda)\Big) t_v \\
&= \Big(\prod_{\lambda \in E}(t_v - t_\lambda t^*_\lambda)\Big) \Big(\prod_{\mu \in \vee (E' H)} (t_v - t_\mu t^*_\mu) + \sum_{\mu \in \vee (E' H)} Q(t)^{\vee (E' H)}_\mu \Big).
\end{align*}
Hence \eqref{eq:orth to FH} gives $\prod_{\lambda \in E}(t_v - t_\lambda t^*_\lambda) = \big(\prod_{\lambda \in E}(t_v - t_\lambda t^*_\lambda)\big) \big(\sum_{\mu \in \vee (E' H)} Q(t)^{\vee (E' H)}_\mu\big)$, and hence belongs to $I_H$ because $\vee (E' H) \subset \L H$, so each $Q(t)^{\vee (E' H)}_\mu \in I_H$.
\end{proof}

Finally, before proving Theorem~\ref{thm:quotient conjecture}, we need to recall some notation and definitions from \cite{RSY2} and \cite{Si1}.

Let $(\L,d)$ be a finitely aligned $k$-graph, and let $G \subset \L$. As in \cite[Definition~3.3]{RSY2}, $\Pi G$ denotes the smallest subset of $\L$ which contains $G$ and has the property that if $\lambda,\mu$ and $\sigma$ belong to $G$ with $d(\lambda) = d(\mu)$ and $s(\lambda) = s(\mu)$ and if $(\alpha,\beta) \in \Lmin(\mu,\sigma)$, then $\lambda\alpha \in G$. If follows from \cite[Lemma~3.2]{RSY2} that $\Pi G$ is finite when $G$ is. We denote by $\Pi G \times_{d,s} \Pi G$ the set of pairs $\{(\lambda,\mu) \in \Pi G \times \Pi G : d(\lambda) = d(\mu), s(\lambda) = s(\mu)\}$.

Let $\{t_\lambda : \lambda \in \L\}$ satisfy (TCK1)--(TCK3). As in \cite[Proposition~3.5]{RSY2}, for a finite set $G \subset \L$ and a path $\lambda \in \Pi G$, we write $Q(t)^{\Pi G}_\lambda$ for the projection 
\begin{equation}\label{eq:Qdef}
Q(t)^{\Pi G}_\lambda := 
\prod_{\lambda\lambda' \in (\Pi G)\setminus \{\lambda\}} (t_\lambda t^*_\lambda - t_{\lambda\lambda'} t^*_{\lambda\lambda'}),
\end{equation}
and for $(\lambda,\mu) \in \Pi G \times_{d,s} \Pi G$, we define
\[
\Theta(t)^{\Pi G}_{\lambda,\mu} := t_\lambda \Big(\prod_{\lambda\lambda' \in (\Pi G)\setminus \{\lambda\}} (t_{s(\lambda)} - t_{\lambda'} t^*_{\lambda'})\Big) t^*_\mu.
\]
By \cite[Lemma~3.10]{RSY2}, we have
\[
Q(t)^{\Pi G}_\lambda t_\lambda t^*_\mu = \Theta(t)^{\Pi G}_{\lambda,\mu} = t_\lambda t^*_\mu Q(t)^{\Pi G}_\mu.
\]

Finally, recall from \cite[Definition~4.4]{Si1} that a graph morphism $x: \Omega_{k,m} \to \L$ is a \emph{boundary path of $\L$} if, whenever $n \le m$ and $E \in x(n)\FE(\L)$, we have $x(n, n+d(\lambda)) = \lambda$ for some $\lambda \in E$. We write $r(x)$ for $x(0)$ and $d(x)$ for $m$. The collection $\partial\L  := \{x : x\text{ is a boundary path of $\L$}\}$ is called the boundary-path space of $\L$. For $\lambda \in \L$ and $x \in \partial\L $ with $r(x) = s(\lambda)$, there is a unique boundary path $\lambda x$ such that $(\lambda x)(0, d(\lambda)) = \lambda$ and $(\lambda x)(d(\lambda), d(\lambda)+n) = x(0,n)$ for all $n \in \NN^k$. Likewise, given $x \in \partial\L $ and $n \le d(x)$, there is a unique boundary path $x|^{d(x)}_n$ such that $(x|^{d(x)}_n)(0, m) = x(n, n+m)$ for all $m \in \NN^k$. As in \cite[Definition~4.6]{Si1}, we define partial isometries $\{S_\lambda : \lambda \in \L\} \subset \Bb(\ell^2(\partial\L ))$ by 
\[
S_\lambda e_x := \delta_{s(\lambda), r(x)} e_{\lambda }.
\] 
Lemma~4.7 of \cite{Si1} shows that $\{S_\lambda : \lambda \in \L\}$ is a Cuntz-Krieger $\L$-family called the \emph{boundary-path representation} and that
\begin{equation}\label{eq:S_lambda*}
S^*_\lambda e_x =
\begin{cases}
e_{x|^{d(x)}_{d(\lambda)}} &\text{ if $x(0, d(\lambda)) = \lambda$} \\
0 &\text{ otherwise.}
\end{cases}
\end{equation}

\begin{proof}[Proof of Theorem~\ref{thm:quotient conjecture}] 
Fix $v \in \L^0 \setminus \L H$ and fix $E \in \FE(\L\setminus\L H)\setminus\Ee_H$.

\begin{claim}\label{clm:projection not in ideal}
Claim~1: For all $a \in \lsp\{s_\lambda s^*_\mu : \lambda,\mu \in \L H\}$, we have
\begin{itemize}
\item[(1)] $\|s_v - a\| \ge 1$; and
\item[(2)] $\big\|\big(\prod_{\lambda \in E} (s_{r(E)} - s_\lambda s^*_\lambda)\big) - a\big\| \ge 1$
\vadjust{\kern-3pt}.
\end{itemize}
\end{claim}
\begin{proof}[\rm Proof of Claim~\ref{clm:projection not in ideal}]
Express $a = \sum_{\lambda \in F} a_{\lambda,\mu} s_\lambda s^*_\mu$ where $F$ is a finite subset of $\L H$, and
$\{a_{\lambda,\mu} : \lambda,\mu \in F\} \subset \CC$. Let $\pi_S$ be the boundary-path representation of $C^*(\L)$ and let $A := \pi_S(a) = \sum_{\lambda,\mu \in F} a_{\lambda,\mu} S_\lambda S^*_\mu$.

To check (1), note that since $v \not\in H$ and since $H$ is saturated, we have that $v F \cap \L^0 = \emptyset$ and 
that $v F \not\in \FE(\L)$. Hence there exists $\tau \in v\L$ such that $\Lmin(\tau,\lambda) = \emptyset$ for all $\lambda \in F$. By \cite[Lemma~4.7(1)]{Si1}, there exists a boundary path $x$ in $s(\tau)\partial\L$. By choice of $\tau$, we have that $\tau x \in v\partial\L \setminus F\partial\L$. But now
\begin{equation}\label{eq:norm is 1}\textstyle
\|S_v - A\| 
\ge \|(S_v - A) e_{\tau x}\| 
= \|S_v e_{\tau x} - \sum_{\lambda,\mu \in F} (a_{\lambda,\mu} S_\lambda S^*_\mu e_{\tau x})\|.
\end{equation}
Since $\tau x \not\in F\partial\L$ by choice, \eqref{eq:S_lambda*} gives $S^*_\mu e_{\tau x} = 0$ for all $\mu \in F$, and hence~\eqref{eq:norm is 1} gives $\|S_v - A\| \ge \|S_v e_{\tau x}\| = \|e_{\tau x}\| = 1$. Since $\pi_S$ is a $C^*$-homomorphism, and hence norm-decreasing, this establishes~(1).

For (2), note that $E \not\in \Ee_H$, and $F \subset \L H$ is finite, so we know that $E \cup F \not\in \FE(\L)$.  Hence there exists $\tau \in \L$ such that $\Lmin(\sigma,\tau) = \emptyset$ for all $\sigma \in E \cup F$.  By \cite[Lemma~4.7(1)]{Si1}, there exists  $x \in \partial\L$ such that $r(x) = s(\tau)$.  Set $y := \tau x \in \partial\L$.  By choice of $\tau$, we have that $y(0, d(\sigma)) \not= \sigma$ for all $\sigma \in E \cup F$.  Hence $S^*_\sigma e_y = 0$ for all $\sigma \in E \cup G$ by~\eqref{eq:S_lambda*}. In particular, $\sigma \in F$ implies $S_\sigma^* e_y = 0$, so $A e_y = 0$, and $\lambda \in E$ implies $S_\lambda^* e_y = 0$. It  follows that $\big(\prod_{\lambda \in E}(S_{r(E)} - S_\lambda S^*_\lambda)\big) e_y = S_{r(E)} e_y = e_y$. Hence
\[
\textstyle 
\big\|\big(\prod_{\lambda \in E}(S_{r(E)} - S_\lambda S^*_\lambda) - A\big)\big\| \ge
\big\|\big(\prod_{\lambda \in E}(S_{r(E)} - S_\lambda S^*_\lambda) - A\big) e_y\big\| = \|e_y\| = 1.
\]
It follows that $\big\|\prod_{\lambda \in E}(S_{r(E)} - S_\lambda S^*_\lambda) - A\big\| \ge 1$. Again since $\pi_S$ 
is norm-decreasing, this establishes~(2).
\hfil\penalty100\hbox{}\nobreak\hfill\hbox{\qed\ Claim~\ref{clm:projection not in ideal}}
\renewcommand\qed{}
\end{proof} 
Since $I_H \subset C^*(\L)$ is fixed under the gauge action, $\gamma$ descends to a strongly continuous action $\theta$ of $\TT^k$ on $C^*(\L) / I_H$ such that $\theta_z \circ \pi^{\Ee_H}_{s + I_H} = \pi^{\Ee_H}_{s + I_H} \circ \gamma_z$ fo all $z \in \TT^k$. 

It is easy to check using~(TCK3) that $\lsp\{s_\lambda s^*_\mu : \lambda,\mu \in \L H\}$ is a dense subset of $I_H$. Hence Claim~\ref{clm:projection not in ideal} shows that neither $s_v$ nor $\prod_{\lambda \in E}(s_{r(E)} - s_\lambda s^*_\lambda)$ belongs to $I_H$. Since $v \in \L^0 \setminus H$ and $E \in \FE(\L\setminus\L H)\setminus \Ee_H$ were arbitrary, and since Lemma~\ref{lem:EsubH its own bar} shows that $\Ee_H$ is satiated, the gauge-invariant uniqueness theorem \cite[Theorem~6.1]{Si1} shows that $\pi^{\Ee_H}_{s + I_H}$ is injective.
\end{proof}

\section{Gauge-invariant ideals in $C^*(\L)$} \label{sec:ideal listing}
Theorem~\ref{thm:quotient conjecture} and \cite[Theorem~6.1]{Si1} combine to show that every nontrivial gauge-invariant ideal in $C^*(\L\setminus\L H;\Ee_H)$ which contains no vertex projection $s_{\Ee_H}(v)$ must contain some collection of projections
\[
\textstyle\big\{\prod_{\lambda \in E} \big(s_{\Ee_H}(r(E)) - s_{\Ee_H}(\lambda)  s_{\Ee_H}(\lambda)^*\big) : E \in B\big\}
\] 
where $B$ is a subset of $\FE(\L\setminus \L H) \setminus \Ee_H$.  

Since $C^*(\L\setminus\L H; \Ee_H)$ itself is the quotient of $C^*(\L)$ by $I_H$, it follows that the ideals $I$ of $C^*(\L)$ such that the set $H_I$ defined in Lemma~\ref{lem:H_I sat,hered} is equal to $H$ should be indexed by some collection of subsets of $\FE(\L\setminus\L H) \setminus \Ee_H$.

In this section, we show that the gauge-invariant ideals of $C^*(\L)$ are indexed by pairs $(H,B)$ where $H$ is a saturated hereditary subset of $\L^0$ and $B$ is a subset of $\FE(\L\setminus \L H) \setminus \Ee_H$ such that $B \cup \Ee_H$ is satiated.

\begin{dfn} 
Let $(\L,d)$ be a finitely aligned $k$-graph and let $H \subset \L^0$ be saturated and hereditary.  Let $B$ be a subset of $\FE(\L\setminus \L H)$.  We define $J_{H,B}$ to be the ideal of $C^*(\L)$ generated by
\[
\textstyle 
\big\{s_v : v \in H\big\} \cup \big\{ \prod_{\lambda \in E} (s_{r(E)} - s_\lambda s^*_\lambda) : E \in B\big\}.
\] 
We define $I(\L\setminus\L H)_B$ to be the ideal of $C^*(\L \setminus \L H; \Ee_H)$ generated by
\[
\textstyle
\big\{\prod_{\lambda \in E} (s_{\Ee_H}(r(E)) - s_{\Ee_H}(\lambda) s_{\Ee_H}(\lambda)^*) : E \in B\big\}.
\]
\end{dfn}

If $H \subset \L^0$ is saturated and hereditary, and if $B$ is a subset of $\FE(\L \setminus\L H)\setminus\Ee_H$ such that $\Ee_H \cup B$ is satiated, then $q(J_{H,B}) \cong I(\L\setminus\L H)_B$ where $q$ is the quotient map from $C^*(\L)$ to $C^*(\L) / I_H \cong C^*(\L\setminus\L H; \Ee_H)$.

We now investigate the structure of $C^*(\L)/J_{H,B}$.

\begin{lemma}\label{lem:quotients equal} 
Let $(\L,d)$ be a finitely aligned $k$-graph and let $H \subset \L^0$ be saturated and hereditary.  Let $B$ be a subset of $\FE(\L\setminus \L H) \setminus \Ee_H$ such that $\Ee_H \cup B$ is satiated.  Then
\[
C^*(\L\setminus \L H; \Ee_H) / I(\L\setminus\L H)_B = C^*(\L\setminus \L H; (\Ee_H \cup B)).
\]
\end{lemma}
\begin{proof} 
By Lemma~\ref{lem:univ-quotients}, we have that $C^*(\L\setminus \L H; \Ee_H) \cong \Tt C^*(\L \setminus\L H) / J_{\Ee_H}$ and $C^*(\L\setminus \L H; (\Ee_H \cup B)) \cong \Tt C^*(\L \setminus\L H) / J_{\Ee_H \cup B}$. Hence we just need to show that $a \in \Tt C^*(\L \setminus\L H)$ belongs to $J_{\Ee_H \cup B}$ if and only if $q(a) \in I(\L\setminus\L H)_B$ where $q : \Tt C^*(\L \setminus\L H) \to C^*(\L\setminus \L H; \Ee_H)$ is the quotient map.

By definition of $I(\L\setminus\L H)_B$, the inverse image $q^{-1}(I(\L\setminus\L H)_B)$ under the quotient map is precisely the ideal in $\Tt C^*(\L \setminus \L H)$ generated by
\[
\begin{split}\textstyle
\{\prod_{\lambda \in E} &(s_{\Tt}(r(E)) - s_\Tt(\lambda) s_\Tt(\lambda)^*) : E \in B\} \\
&\textstyle\cup\ \{\prod_{\lambda \in E} (s_\Tt(r(E)) - s_\Tt(\lambda) s_\Tt(\lambda)^*) : E \in \Ee_H\};
\end{split}
\]
that is, $q^{-1}(I(\L\setminus\L H)_B) = J_{\Ee_H \cup B}$ as required.
\end{proof}

\begin{cor} \label{cor:2nd quotient a relCK} 
Let $(\L,d)$ be a finitely aligned $k$-graph, let $H \subset \L^0$ be saturated and hereditary, and let $B \subset \FE(\L\setminus \L H) \setminus \Ee_H$.  Then
\[
C^*(\L)/J_{H,B} \cong C^*(\L\setminus \L H; (\Ee_H \cup B)).
\]
\end{cor}
\begin{proof} 
We will show that $C^*(\L)/ J_{H,B} = (C^*(\L)/I_H)/I(\L\setminus\L H)_B$; the result then follows from Lemma~\ref{lem:quotients equal}.  Let
\begin{align*}
q_{H,B} &: C^*(\L) \to C^*(\L)/ J_{H,B}, \\
q_H &: C^*(\L) \to C^*(\L) / I_H, \\
q_B &: C^*(\L) / I_H \to (C^*(\L) / I_H)/ I(\L\setminus\L H)_B 
\end{align*} 
be the quotient maps.  It is clear that the kernel of $q_{H,B}$ is contained in that of $q_B \circ q_H$, giving a canonical homomorphism $\pi_1$ of $C^*(\L)/ J_{H,B}$ onto $(C^*(\L)/I_H)/I(\L\setminus\L H)_B$.  On the other hand, since $I_H \subset J_{H,B}$, there is a canonical homomorphism $\pi_2$ of $C^*(\L)/I_H$ onto $C^*(\L)/ J_{H,B}$ whose kernel contains $I(\L\setminus\L H)_B$ by definition.  It follows that $\pi_2$ descends to a canonical homomorphism $\tilde\pi_2$ of $(C^*(\L)/I_H)/I(\L\setminus\L H)_B$ onto $C^*(\L)/ J_{H,B}$ which is inverse to $\pi_1$.
\end{proof}

\begin{dfn} 
Let $(\L,d)$ be a finitely aligned $k$-graph.  For each gauge-invariant ideal $I$ in $C^*(\L)$, recall that $H_I$ denotes $\{v \in \L^0 : s_v \in I\}$, and define
\[
\textstyle B_I := \big\{E \in \FE(\L \setminus \L H_I) \setminus \Ee_{H_I} : 
\prod_{\lambda \in E}  (s_{\Ee_{H_I}}(r(E)) - s_{\Ee_{H_I}}(\lambda) s_{\Ee_{H_I}}(\lambda)^*) \in q_{H_I}(I)\big\},
\]
where $q_{H_I}$ is the quotient map from $C^*(\L)$ to $C^*(\L) / I_{H_I}$.
\end{dfn}

\begin{theorem}\label{thm:every g-i ideal is a JH,B} 
Let $(\L,d)$ be a finitely aligned $k$-graph.
\begin{itemize}
\item[(1)] Let $I$ be a gauge-invariant ideal of $C^*(\L)$.  Then $H_I \subset \L^0$ is nonempty saturated and hereditary, $\Ee_{H_I} \cup B_I$ is a satiated subset of $\FE(\L \setminus \L {H_I})$, and $I = J_{H_I, B_I}$.  
\item[(2)] Let $H \subset \L^0$ be nonempty, saturated and hereditary, and let $B$ be a subset of $\FE(\L\setminus \L H) \setminus \Ee_H$ such that $\Ee_H \cup B$ is satiated in $\L \setminus \L H$.  Then $H_{J_{H,B}} = H$ and $B_{J_{H,B}} = B$.
\end{itemize}
\end{theorem}
\begin{proof} Theorem~6.1 of \cite{Si1} shows that $H_I$ is nonempty, and Lemma~\ref{lem:H_I sat,hered} shows that it is saturated and hereditary.  That $\Ee_H \cup B_I$ is satiated follows from \cite[Corollary~4.10]{Si1}.

Let $I$ be a gauge-invariant ideal of $C^*(\L)$. We have $J_{H_I, B_I} \subset I$ by definition, so there is a canonical homomorphism $\pi$ of $C^*(\L) / J_{H_I, B_I}$ onto $C^*(\L)/I$.  By Corollary~\ref{cor:2nd quotient a relCK}, this gives us a homomorphism, also denoted $\pi$ of $C^*(\L\setminus\L H_I; \Ee_{H_I} \cup B_I)$ onto $C^*(\L)/I$.  Since $I$ is gauge-invariant, the gauge action on $C^*(\L)$ descends to an action $\theta$ of $\TT^k$ on $C^*(\L)/I$ such that $\theta_z \circ \pi = \pi \circ \gamma_z$ where $\gamma$ is the gauge action on $C^*(\L\setminus\L H_I; \Ee_{H_I} \cup B_I)$.

Suppose that $\pi(s_{\Ee_{H_I} \cup B_I}(v))$ is equal to $0$ in $C^*(\L)/I$. Then $s_v \in I$ by definition, so $v \in H_I$.  Hence $\pi(s_{\Ee_{H_I} \cup B_I}(v)) \not= 0$ for all $v \in (\L \setminus \L H_I)^0$.

Now suppose that $E \in \FE(\L \setminus \L H_I)$ satisfies 
\[
\textstyle
\pi\Big(\prod_{\lambda \in E} (s_{\Ee_{H_I} \cup B_I}(r(E)) - s_{\Ee_{H_I} \cup B_I}(\lambda) 
s_{\Ee_{H_I} \cup B_I}(\lambda)^*)\Big) = 0_{C^*(\L)/I}.
\]
Then either $E \in \Ee_{H_I}$, or else $E \in B_I$ by the definition of $B_I$.  But then $\prod_{\lambda \in E} (s_{r(E)} - s_\lambda s^*_\lambda) \in J_{H_I, B_I}$, so that
\[
\prod_{\lambda \in E} (s_{\Ee_{H_I} \cup B_I}(r(E)) - s_{\Ee_{H_I} \cup B_I}(\lambda) s_{\Ee_{H_I} \cup B_I}(\lambda)^*) = 0_{C^*(\L\setminus \L H_I; \Ee_{H_I} \cup B_I)}.
\] 
Hence $\pi\Big(\prod_{\lambda \in E} (s_{\Ee_{H_I} \cup B_I}(r(E)) - s_{\Ee_{H_I} \cup B_I}(\lambda) s_{\Ee_{H_I} \cup B_I}(\lambda)^*)\Big) \not= 0$ for all $E \in \FE(\L) \setminus (\Ee_H \cup B)$.

By the previous three paragraphs we can apply \cite[Theorem~6.1]{Si1} to see that $\pi$ is faithful, and hence that $I = J_{H_I, B_I}$ as required.

Now let $H \subset \L^0$ be saturated and hereditary, and let $B$ be a subset of $\FE(\L \setminus\L H)\setminus \Ee_H$ such that $\Ee_H \cup B$ is satiated.

We have $H \subset H_{J_{H,B}}$ and $B \subset B_{J_{H,B}}$ by definition.  If $v \in H_{J_{H,B}}$, then $s_v \in J_{H,B}$ and hence its image in $C^*(\L \setminus \L H; \Ee_H \cup B)$ is trivial.  It follows that either $v \in H$ or $s_{\Ee_H \cup B}(v) = 0$.  But $s_{\Ee_H \cup B}(v) \not= 0$ for all $v \in (\L \setminus \L H)^0$ by \cite[Theorem~4.3]{Si1}, giving $v \in H$.

If $E \in B_{J_{H,B}}$, then we have
\[
\prod_{\lambda \in E} (s_{\Ee_H}(v) - s_{\Ee_H}(\lambda) s_{\Ee_H}(\lambda)^*) \in I(\L\setminus\L H)_B \subset C^*(\L\setminus\L H; \Ee_H).
\]
Hence $\prod_{\lambda \in E} (s_{\Ee_H \cup B}(v) - s_{\Ee_H \cup B}(\lambda) s_{\Ee_H \cup B}(\lambda)^*)$ is equal to the zero element of $C^*(\L\setminus\L H; \Ee_H) / I(\L\setminus\L H)_B = C^*(\L\setminus\L H; \Ee_H \cup B)$.  Since $\Ee_H \cup B$ is satiated, it follows that either $E \in \Ee_H$ or $E \in B$ by \cite[Theorem~4.3]{Si1}.  But $B_{J_{H,B}} \cap \Ee_H = \emptyset$ by definition, and it follows that $E \in B$ as required.
\end{proof}

\begin{rmk} \begin{itemize}
\item[(1)] Given a saturated hereditary $H \subset \L^0$, the ideal $I_H$ (see Notation~\ref{ntn:I_H}) is listed by Theorem~\ref{thm:every g-i ideal is a JH,B} as $J_{H,\emptyset}$.

\item[(2)] It seems difficult to establish an analogue of Lemma~\ref{lem:Lambda(H) Vs I_H} for arbitrary $J_{H,B}$. A good strategy would be to aim to describe $I(\L\setminus\L H)_B = J_{H,B} / I_H$ as (Morita equivalent to) a $k$-graph algebra. But this seems difficult even when $B$ is ``singly generated:'' i.e. when $\Ee_H \cup B$ is the satiation (see \cite[Definition~5.1]{Si1}) of $\Ee_H \cup \{E\}$ where $E \in \FE(\L\setminus \L H) \setminus \Ee_H$.
\end{itemize}\end{rmk}

\section{The lattice order}\label{sec:lattice}

In this section we describe the lattice ordering of the gauge-invariant ideals of $C^*(\L)$ in terms of a lattice order on the pairs $(H,B)$ where $H \subset \L^0$ is saturated and hereditary, and $B$ is a subset of $\FE(\L\setminus \L H) \setminus \Ee_H$ such that $\Ee_H \cup B$ is satiated.

\begin{dfn} \label{dfn:(H,B) order}
Let $(\L,d)$ be a finitely aligned $k$-graph. Define 
\[\begin{split}
\Ipairs :=
\big\{(H,B) : {}&\emptyset\not= H \subset \L^0, H\text{ is saturated and hereditary, }\\
&B \subset \FE(\L\setminus\L H)\setminus\Ee_H\text{ and }\Ee_H \cup B\text{ is satiated}\big\}.
\end{split}\]
Define a relation $\preceq$ on $\Ipairs$ by $(H_1, B_1) \preceq (H_2, B_2)$ if and only if 
\begin{itemize}
\item[(1)] $H_1 \subset H_2$; and 
\item[(2)] if $E \in B_1$ and $r(E) \not\in H_2$, then $E \setminus EH_2$ belongs to $\Ee_{H_2} \cup B_2$.
\end{itemize}
\end{dfn}

\begin{theorem} \label{thm:gauge-invariant ideals} 
Let $(\L,d)$ be a finitely aligned $k$-graph. The map $(H,B) \mapsto J_{H,B}$ is a lattice isomorphism between $(\Ipairs, \preceq)$ and $(I^\gamma(\L), \subset)$ where $I^\gamma(\L)$ denotes the collection of gauge-invariant ideals of $C^*(\L)$. 
\end{theorem}
\begin{proof} 
Theorem~\ref{thm:every g-i ideal is a JH,B} implies that $(H,B) \mapsto J_{H,B}$ is a bijection between $\Ipairs$ and $I^\gamma(C^*(\L))$. Hence, we need only establish that for $(H_1, B_1), (H_2, B_2) \in \Ipairs$,
\begin{equation}\label{eq:subset vs preceq}
\text{$J_{H_1, B_1} \subset J_{H_2, B_2}$ if and only if $(H_1, B_1) \preceq (H_2, B_2)$.}
\end{equation}

First suppose that $J_{H_1, B_1} \subset J_{H_2, B_2}$. Theorem~\ref{thm:every g-i ideal is a JH,B} shows immediately that $H_1 \subset H_2$, so if we can show that $F \in B_1$ with $r(F) \not\in H_2$ implies $F\setminus FH_2 \in \Ee_{H_2} \cup B_2$, it will follow that $(H_1, B_1) \preceq (H_2, B_2)$. 

Suppose that $E = F\setminus FH_2$ for some $F \in B_1$ with $r(F) \not\in H_2$. Suppose further for contradiction that $E \not\in \Ee_{H_2} \cup B_2$. Let $q_i : C^*(\L) \to C^*(\L)/ J_{H_i, B_i}$ where $i = 1,2$ denote the quotient maps; by Corollary~\ref{cor:2nd quotient a relCK}, we can regard $q_i$ as a homomorphism of $C^*(\L)$ onto $C^*(\L\setminus\L H_i; \Ee_{H_i} \cup B_i)$ for $i=1,2$. Since $J_{H_1, B_1} \subset J_{H_2, B_2}$, there is a homomorphism $\pi : C^*(\L\setminus \L H_1; \Ee_{H_1}\cup B_1) \to C^*(\L\setminus \L H_2; \Ee_{H_2}\cup B_2)$ such that $\pi\circ q_1 = q_2$. Since $F \in B_1$, we have $q_1\big(\prod_{\lambda \in F} (s_{r(F)} - s_\lambda s^*_\lambda)\big) = 0$, and hence
\begin{equation}\label{eq:q2 image 0}\textstyle
q_2\big(\prod_{\lambda \in F} (s_{r(F)} - s_\lambda s^*_\lambda)\big) = \pi\Big(q_1\big(\prod_{\lambda \in F} (s_{r(F)} - s_\lambda s^*_\lambda)\big)\Big) = 0.
\end{equation} 
Since $s(\lambda) \in H_2$ implies $q_2(s_\lambda s^*_\lambda) = 0$ by definition, we have that 
\begin{equation}\label{eq:q2 image}\textstyle
q_2\Big(\prod_{\lambda \in F} (s_{r(F)} - s_\lambda s^*_\lambda)\Big)
= \prod_{\lambda \in E}\big(s_{\Ee_{H_2} \cup B_2}(r(E)) - s_{\Ee_{H_2} \cup B_2}(\lambda)s_{\Ee_{H_2} \cup
B_2}(\lambda)^*\big),
\end{equation}

We consider two cases:

Case~1: $E$ belongs to $\FE(\L\setminus\L H_2)$. Then since $E \not\in \Ee_{H_2} \cup B_2$, \cite[Corollary~4.10]{Si1} ensures that $\prod_{\lambda \in E}\big(s_{\Ee_{H_2} \cup B_2}(r(E)) - s_{\Ee_{H_2} \cup B_2}(\lambda) s_{\Ee_{H_2} \cup B_2}(\lambda)^*\big)$ is nonzero.

Case~2: $E \not\in \FE(\L\setminus\L H_2)$. Then there exists $\mu \in r(E)\L \setminus \L H_2$ with $\Ext(\mu;E) = \emptyset$; we then have 
\[\begin{split}\textstyle
\prod_{\lambda \in E}\big(s_{\Ee_{H_2} \cup B_2}(r(E)) - s_{\Ee_{H_2} \cup B_2}(\lambda)s_{\Ee_{H_2} \cup
B_2}(\lambda)^*\big) &s_{\Ee_{H_2} \cup B_2}(\mu) s_{\Ee_{H_2} \cup B_2}(\mu)^* \\
&= s_{\Ee_{H_2} \cup B_2}(\mu) s_{\Ee_{H_2} \cup B_2}(\mu)^*
\end{split}\]
by (TCK3). Since $s_{\Ee_{H_2} \cup B_2}(\mu) s_{\Ee_{H_2} \cup B_2}(\mu)^* \not= 0$ by \cite[Corollary~4.10]{Si1}, it follows that
\[
\prod_{\lambda \in E}\big(s_{\Ee_{H_2} \cup B_2}(r(E)) - s_{\Ee_{H_2} \cup B_2}(\lambda)s_{\Ee_{H_2} \cup
B_2}(\lambda)^*\big) s_{\Ee_{H_2} \cup B_2}(\mu) s_{\Ee_{H_2} \cup B_2}(\mu)^* \not= 0.
\]

In either case,~\eqref{eq:q2 image} shows that $q_2\big(\prod_{\lambda \in F} (s_{r(F)} - s_\lambda s^*_\lambda)\big)$ is nonzero, contradicting~\eqref{eq:q2 image 0}. This establishes the ``only if'' assertion of~\eqref{eq:subset vs preceq}.

Now suppose that $(H_1, B_1) \preceq (H_2, B_2) \in \Ipairs$. Let $v \in H_1$. Since $(H_1, B_1) \preceq (H_2, B_2)$, we have that $H_1 \subset H_2$, and hence $v \in H_2$ giving $s_v \in J_{H_2, B_2}$ by definition. Now let $E \in B_1$. If $r(E) \in H_2$, then $s_{r(E)} \in J_{H_2, B_2}$ by definition, and hence $\prod_{\lambda \in E} (s_{r(E)} - s_\lambda s^*_\lambda) = \big(\prod_{\lambda \in E} (s_{r(E)} - s_\lambda s^*_\lambda)\big) s_{r(E)} \in J_{H_2, B_2}$. If $r(E) \not\in H_2$, then since $(H_1, B_1) \preceq (H_2, B_2)$, we have that $E \setminus E H_2 \in
\Ee_{H_2} \cup B_2$. For $\lambda \in \L H_2$, we have $s_\lambda s^*_\lambda  = s_\lambda s_{s(\lambda)} s^*_\lambda \in J_{H_2, B_2}$ and hence $q_2(s_\lambda s^*_\lambda) = 0$, so
\begin{equation}
\label{eq:quotient image}
q_2\Big(\prod_{\lambda \in E} (s_{r(E)} - s_\lambda s^*_\lambda)\Big) = \prod_{\lambda \in E\setminus E H_2}
(s_{\Ee_{H_2}\cup B_2}(r(E)) - s_{\Ee_{H_2}\cup B_2}(\lambda) s_{\Ee_{H_2}\cup B_2}(\lambda)^*).
\end{equation}
Since $E \setminus E H_2 \in \Ee_{H_2} \cup B_2$, and since $\{s_{\Ee_{H_2} \cup B_2}(\lambda) : \lambda \in \L \setminus \L H_2\}$ is a relative Cuntz-Krieger $(\L\setminus\L H_2; E_{H_2} \cup B_2)$-family, relation~(CK) gives 
\[\textstyle
\prod_{\lambda \in E \setminus E H_2} (s_{\Ee_{H_2} \cup B_2}(r(E)) - s_{\Ee_{H_2} \cup B_2}(\lambda) s_{\Ee_{H_2} \cup B_2}(\lambda)^*) = 0.
\]
Hence $\prod_{\lambda \in E} (s_{r(E)} - s_\lambda s^*_\lambda) \in \ker q_2 = J_{H_2, B_2}$ by~\eqref{eq:quotient image} and Corollary~\ref{cor:2nd quotient a relCK}. 

Since all the generating projections of $J_{H_1, B_1}$ belong to $J_{H_2, B_2}$, it follows that $J_{H_1, B_1} \subset J_{H_2, B_2}$, establishing the ``if'' assertion of~\eqref{eq:subset vs preceq}.
\end{proof}

\section{$k$-graphs in which all ideals are gauge-invariant} \label{sec:CK ideals}
In this section we use the Cuntz-Krieger uniqueness theorem of \cite{Si1} to show that for a certain class of $k$-graphs, the ideals $J_{H,B}$ identified in Section~\ref{sec:ideal listing} are all the ideals in $C^*(\L)$; that is, every ideal in $C^*(\L)$ is gauge-invariant.

Recall from \cite[Definition~6.2]{Si1} that if $x : \Omega_{k, d(x)} \to \L$ and $y : \Omega_{k, d(y)} \to \L$ are graph morphisms, then $\MCE(x,y)$ is the collection of all graph morphisms $z : \Omega_{k, d(z)} \to \L$ such that $d(z)_i = \max{d(x)_i, d(y)_i}$ for $1 \le i \le k$, and such that $z|_{\Omega_{k, d(x)}} = x$ and $z|_{\Omega_{k, d(y)}} = y$.

Recall also from \cite[Theorem~6.3]{Si1} that if $(\L,d)$ is a finitely aligned $k$-graph and $\Ee$ is a subset of $\FE(\L)$, then $(\L,\Ee)$ is said to satisfy \emph{condition~{\rm(C)}} if
\begin{itemize}
\item[(1)] For all $v \in \L^0$ there exists $x \in v\partial(\L;\Ee)$ such that for distinct $\lambda,\mu$ in $\L r(x)$, we have $\MCE(\lambda x, \mu x) = \emptyset$; and
\item[(2)] for each $F \in v\FE(\L)\setminus\overline\Ee$, there is a path $x$ as in~(1) such that $x \in v\partial(\L;\Ee) \setminus F\partial(\L;\Ee)$. 
\end{itemize}

\begin{dfn}
Let $(\L,d)$ be a finitely aligned $k$-graph. We say that $\L$ satisfies \emph{condition~\eqref{eq:new condition K}} if 
\begin{equation}\label{eq:new condition K}
\text{$(\L \setminus \L H, \Ee_H)$ satisfies condition~(C) for each saturated, hereditary $H \subset \L^0$.}\tag{D}
\end{equation}
\end{dfn}

\begin{theorem}\label{thm:cond(D)--> every ideal a JH,B} 
Let $(\L,d)$ be a finitely aligned $k$-graph which satisfies condition~\eqref{eq:new condition K}.
\begin{itemize}
\item[(1)] Let $I$ be an ideal of $C^*(\L)$.  Then $H_I$ is nonempty, saturated and hereditary, $B_I \cup \Ee_{H_I}$ is satiated in $\L \setminus \L {H_I}$, and $I = J_{H_I, B_I}$.  
\item[(2)] Let $H \subset \L^0$ be nonempty, saturated and hereditary, and let $B \subset \FE(\L\setminus \L H) \setminus \Ee_H$ be such that $B \cup \Ee_H$ is satiated in $\L \setminus \L H$.  Then $H_{J_{H,B}} = H$ and $B_{J_{H,B}} = B$.
\end{itemize}
\end{theorem}
\begin{proof}
The proof of~(1) is the same as the proof of of Theorem~\ref{thm:every g-i ideal is a JH,B}(1) except that, since we do not know \emph{a priori} that $I$ is gauge-invariant, we do not automatically have an action $\pi$ on $C^*(\L)/I$ such that $\theta_z \circ \pi = \pi \circ \gamma_z$. Consequently, we cannot apply \cite[Theorem~6.1]{Si1} to deduce that $\pi$ is faithful; instead, we use our assumption that $(\L\setminus\L H, \Ee_H)$ satisfies condition~(C) to apply \cite[Theorem~6.3]{Si1}.

The proof of~(2) is identical the to proof of part~(2) of Theorem~\ref{thm:every g-i ideal is a JH,B}.
\end{proof}

\section{Classifiability} \label{ch:classifiable} 
In this section we investigate when $C^*(\L)$ is a Kirchberg-Phillips algebra. We show that all relative $k$-graph algebras $C^*(\L;\Ee)$ fall into the bootstrap class $\mathcal{N}$ of \cite{RSc}.  We show that if $\L$ satisfies condition~(C), then $C^*(\L)$ is simple if and only if $\L$ is \emph{cofinal}. Finally, we show that if in addition every vertex of $\L$ can be reached from a \emph{loop with an entrance}, then $C^*(\L)$ is purely infinite.

The main results in this section are generalisations to arbitrary finitely aligned $k$-graphs of the corresponding results of Kumjian and Pask for row-finite $k$-graphs with no sources in \cite{KP}.

The author would like to thank D. Gwion Evans for drawing his attention to the results of \cite{PQR} which provide the necessary technical machinery for the proof of Proposition~\ref{prp:rel algs nuclear}.

\begin{prop}\label{prp:rel algs nuclear} 
Let $(\L,d)$ be a finitely aligned $k$-graph and let $\Ee$ be a subset of\/ $\FE(\L)$. Then $C^*(\L;\Ee)$ is stably isomorphic to a crossed product of an AF algebra by $\ZZ^k$, and hence falls into the bootstrap class $\mathcal{N}$ of \cite{RSc}; in particular, $C^*(\L;\Ee)$ is nuclear and satisfies the Universal Coefficient Theorem.
\end{prop}

This proposition generalises \cite[Theorem~5.5]{KP}, and the overall strategy of the proof is the same, but the technical details are more complicated, and draw on \cite{RSY2} and \cite{PQR}. We first need to establish some preliminary lemmas, the first of which generalises \cite[Lemma~5.4]{KP}.

\begin{lemma}\label{lem:skew-product AF}
Let $(\L,d)$ be a finitely aligned $k$-graph and let $\Ee \subset \FE(\L)$. Suppose there is a function $b : \L^0 \to \ZZ^k$ such that $d(\lambda) = b(s(\lambda)) - b(r(\lambda))$ for all $\lambda \in \L$. Then $C^*(\L;\Ee)$ is AF.
\end{lemma}
\begin{proof}
The proof is based heavily on that of \cite[Lemma~3.2]{RSY2}.

It suffices to show that for $E \subset \L$ finite, we have that $C^*(\{s_\Ee(\lambda) : \lambda \in E\})$ is finite dimensional. Recalling the definition of $\vee E$ from Notation~\ref{ntn:vee E}, define a map $M$ on finite subsets of $\L$ by 
\begin{equation}\label{eq:Pi tilde}
\begin{split}
M(E) := \{(\lambda_1(0, d(\lambda_1))\lambda_2(n_2, d(\lambda_2)) \dots &\lambda_l(n_l, d(\lambda_l)) : \\
&l \in \NN\setminus\{0\}, \lambda_i \in \vee E, n_i \le d(\lambda_i)\}.
\end{split}
\end{equation}
We claim that
\begin{itemize}
\item[(a)] $M(E)$ is finite;
\item[(b)] $E \subset \vee E \subset M(E)$;
\item[(c)] $\bigvee_{\lambda \in M(E)} b(s(\lambda)) = \bigvee_{\mu \in E} b(s(\mu))$;
\item[(d)] $\lambda,\mu,\sigma,\tau \in E$ implies $s_\Ee(\lambda) s_\Ee(\mu)^* s_\Ee(\sigma) s_\Ee(\tau)^* \in \lsp\{s_\Ee(\eta) s_\Ee(\zeta)^* : \eta,\zeta \in M(E)\}$; and
\item[(e)] if $M^2(E) \not= M(E)$, then $\min\{\sum^k_{i=1}b(s(\lambda))_i : \lambda \in M^2(E) \setminus M(E)\}$ is strictly greater than $\min\{\sum^k_{i=1}b(s(\mu))_i : \mu \in M(E) \setminus E\}$.
\end{itemize}
For~(a), note that each path in $M(E)$ can be factorised as $\alpha_1\dots\alpha_{|d(\lambda)|}$ where each $\alpha_i = \mu(n, n+e_l)$ for some $n \in \NN^k$, $1 \le l  \le k$, and $\mu \in \vee E$. Moreover, $i < j \implies b(s(\alpha_i)) < \big(b(s(\alpha_i)) + d(\alpha_j)\big) \le b(s(\alpha_j)) \implies \alpha_i \not= \alpha_j$. Since $\vee E$ is finite, the number of possible values for $\alpha_i$ is finite, and it follows that $M(E)$ is finite.

We have $E \subset \vee E$ by definition, and $\vee E \subset M(E)$ by taking $l = 1$ in~\eqref{eq:Pi tilde}, establishing~(b).

For~(c), first note that $\lambda \in M(E) \implies s(\lambda) = s(\mu)$ for some $\mu \in \vee E$, so
\begin{equation}\label{eq:first le}
\textstyle \bigvee_{\lambda \in M(E)} b(s(\lambda)) \le \bigvee_{\mu \in \vee E} b(s(\mu)).
\end{equation}
Next recall from \cite[Definition~8.3]{RS1} that for finite $F \subset \L$ 
\[
\MCE(F) := \{\lambda \in \L : d(\lambda) = \bigvee_{\mu \in F} d(\mu), \lambda(0, d(\mu)) = \mu\text{ for all } \mu \in F\},
\]
and that $\vee E = \bigcup\{\MCE(F) : F \subset E\}$. So $\lambda \in \vee E \implies \lambda \in \MCE(F)$ for some subset $F$ of $E$. In particular, $\MCE(F)$ is nonempty, so we must have $F \subset v\L$ for some $v \in \L^0$. Write $n$ for $b(v)$, and calculate:
\[\textstyle
b(s(\lambda)) = n + \bigvee_{\mu \in F} d(\mu) = n + \bigvee_{\mu \in F} (b(s(\mu)) - n) = \bigvee_{\mu \in F} b(s(\mu)).
\]
Hence $\bigvee_{\lambda \in \vee E} b(s(\lambda)) \le \bigvee_{\mu \in E} b(s(\mu))$, so $\bigvee_{\lambda \in M(E)} b(s(\lambda)) \le \bigvee_{\mu \in E} b(s(\mu))$ by~\eqref{eq:first le}. The reverse inequality follows from~(b), establishing~(c).

Claim~(d) follows from~\eqref{eq:Pi tilde} and (TCK3). Finally,~(e) follows from an argument identical to the proof of~(e) in \cite[Lemma~3.2]{RSY2} but with $d(\lambda)$ replaced with $b(\lambda)$ throughout. This establishes the claim.

It now follows as in \cite[Lemma~3.2]{RSY2} that $M^\infty(E) := \bigcup^\infty_{i=1} M^i(E)$ is finite and that $\lsp\{s_\Ee(\lambda) s_\Ee(\mu)^* : \lambda,\mu \in M^\infty(E)\}$ is a finite-dimensional subalgebra of $C^*(\L;\Ee)$ containing $C^*(\{s_\Ee(\lambda) : \lambda \in E\})$.
\end{proof}

Let $\L \times_d \ZZ^k$ be the skew-product $k$-graph which is equal, as a set, to $\L \times \ZZ^k$ and has range, source and degree maps given by $r(\lambda,n) := (r(\lambda), n - d(\lambda))$, $s(\lambda,n) := (s(\lambda), n)$, and $d(\lambda,n) := d(\lambda)$ (see \cite[Definition~5.1]{KP}). For $E \in \Ee$ and $n \in \ZZ^k$, let $E \times_d \{n\} := \{(\lambda, n + d(\lambda)) : \lambda \in E\}$, and let $\Ee \times_d \ZZ^k := \{E \times_d \{n\} : E \in \Ee, n \in \ZZ^k\}$. 

Recall that a \emph{coaction} $\delta$ of a group $G$ on a $C^*$-algebra $A$ is an injective unital homomorphism $\delta : A \to A \otimes C^*(G)$ satisfying the \emph{cocycle identity} $(\id \otimes \delta_G)\circ\delta =  (\delta \otimes \id)\circ\delta$. The \emph{fixed point algebra} is the subspace $A^\delta := \{a \in A : \delta(a) = a \otimes e\}$. There is a universal crossed product algebra $A \times_\delta G$ associated to the triple $(A,G,\delta)$, and this algebra admits a \emph{dual action} $\hat\delta$ of $G$. Crossed product duality says that $A \times_\delta G \times_{\hat\delta} G \cong A \otimes \ell^2(G)$.

The following lemma generalises~\cite[Theorem~7.1]{PQR} to relative $k$-graph algebras.

\begin{lemma} \label{lem:coactions}
Let $(\L,d)$ be a finitely aligned $k$-graph, and let $\Ee$ be a subset of $\FE(\L)$. Then
\begin{itemize}
\item[(1)] $\Ee \times_d \ZZ^k$ is a subset of $\FE(\L \times_d \ZZ^k)$;
\item[(2)] $C^*(\L \times_d \ZZ^k; \Ee \times_d \ZZ^k)$ is AF;
\item[(3)] there is a unique coaction $\delta$ of $\ZZ^k$ on $C^*(\L;\Ee)$ such that $\delta(s_\Ee(\lambda)) := s_\Ee(\lambda) \otimes d(\lambda)$ for all $\lambda \in \L$; and
\item[(4)] the crossed product algebra $C^*(\L;\Ee) \times_\delta \ZZ^k$ is isomorphic to $C^*(\L \times_d \ZZ^k; \Ee \times_d \ZZ^k)$.
\end{itemize}
\end{lemma}
\begin{proof}
For part~(1), fix $E \times_d \{n\} \in \Ee \times_d \ZZ^k$, and suppose that $r(\lambda,m) = r(E \times_d \{n\})$. Then $m = n + d(\lambda)$ and $r(\lambda) = r(E)$. Since $E \in \FE(\L)$, there exists $\alpha \in \Ext(\lambda;E)$. It is straightforward to check that $(\alpha, m + d(\alpha)) \in \Ext((\lambda,m); E \times_d \{n\})$. Since $(\lambda,m)$ was arbitrary, it follows that $E \times_d \{n\} \in \FE(\L \times_d \ZZ^k)$, and since $E \times_d \{n\}$ was itself arbitrary in $\Ee \times_d \ZZ^k$, this establishes~(1).

For~(2), define $b : (\L \times_d \ZZ^k)^0 \to \ZZ^k$ by $b(\lambda,n) := n$. Then the pair $(\L \times_d \ZZ^k, b)$ satisfies the hypotheses of Lemma~\ref{lem:skew-product AF}, so $C^*(\L \times_d \ZZ^k; \Ee \times_d \ZZ^k)$ is AF.

Parts (3)~and~(4) now follow exactly as (i)~and~(ii) of \cite[Theorem~7.1]{PQR}.
\end{proof}

\begin{proof}[Proof of Proposition~\ref{prp:rel algs nuclear}]
We have that $C^*(\L;\Ee) \times_\delta \ZZ^k \cong C^*(\L \times_d \ZZ^k; \Ee \times_d \ZZ^k)$ is AF. But crossed product duality gives $C^*(\L;\Ee) \otimes \ell^2(\ZZ^k) \cong C^*(\L;\Ee) \times_\delta \ZZ^k \times_{\hat\delta} \ZZ^k$, so $C^*(\L;\Ee)$ is stably isomorphic to a crossed product of an AF algebra by $\ZZ^k$.
\end{proof}

Our simplicity result is a direct generalisation of \cite[Proposition~4.8]{KP}, though our proof is based on that of \cite[Proposition~5.1]{BPRS}.

\begin{dfn} 
Let $(\L,d)$ be a finitely aligned $k$-graph.  We say that $\L$ is \emph{cofinal\/} if for all $v \in \L^0$ and $x \in \partial\L $, there exists $n \le d(x)$ such that $v \L x(n) \not= \emptyset$.
\end{dfn}

\begin{prop}\label{prp:simple graph alg}
Let $(\L,d)$ be a finitely aligned $k$-graph, and suppose that $\L$ satisfies condition~{\rm(C)}.  Then $C^*(\L)$ is simple if and only if $\L$ is cofinal.
\end{prop}
\begin{proof} 
First suppose that $\L$ is cofinal, and suppose that $I$ is an ideal in $C^*(\L)$.  If $s_v \in I$ for all $v \in \L^0$, then $I = C^*(\L)$ by (TCK2).  Suppose that $v \in \L^0$ with $s_v \not \in I$.  We must show that $H_I$ is empty, for if so then \cite[Theorem~6.3]{Si1} shows that $I$ is trivial. Since $H_I$ is saturated, we have that 
\begin{equation}\label{eq:saturation consequence}
\text{if $v' \not\in H_I$ and $E \in v\FE(\L)$, then there exists $\lambda \in E$ such that $s(\lambda) \not \in H_I$.}
\end{equation}
To prove the proposition, we first establish the following claim:

\begin{claim}\label{clm:boundary path}
There exists a path $x \in \partial\L$ such that $x(n) \not\in H_I$ for all $n \le d(x)$.
\end{claim}
\begin{proof}[Proof of Claim~\ref{clm:boundary path}]
The proof of the claim is very similar to the proof of \cite[Lemma~4.7(1)]{Si1}, but with minor technical changes needed to establish that we can obtain $x(n) \not\in H_I$ for all $n$. Consequently, we give a proof sketch with frequent references to the proof in \cite{Si1}.

As in the proof of \cite[Lemma~4.7(1)]{Si1}, let $P : \NN^2 \to \NN$ be the position function associated to the diagonal listing of $\NN^2$: 
\[
P(0,0) = 0,\quad P(0,1) = 1,\quad P(1,0) = 2,\quad P(0,2) = 3,\quad P(1,1) = 4,\quad \dots
\]
For $l \in \NN$, let $(i_l, j_l)$ be the unique element of $\NN^2$ such that $P(i_l, j_l) = l$.

We will show by induction that there exists a sequence $\{\lambda_l : l \ge 0\} \subset v\L$ and enumerations $\{E_{l,j} : j \ge 0\}$ of $s(\lambda_l)\FE(\L)$ for all $l \ge 0$ such that
\begin{itemize}
\item[(i)] $s(\lambda_l) \not\in H_I$ for all $l$;
\item[(ii)] $\lambda_{l+1}(0, d(\lambda_l)) = \lambda_l$ for all $l \ge 1$; and
\item[(iii)] $\lambda_{l+1}(d(\lambda_{i_l}), d(\lambda_{l+1})) \in E_{i_l, j_l}\L$ for all $l \ge 0$.
\end{itemize}
As in the proof of \cite[Lemma~4.7(1)]{Si1}, we proceed by induction on $l$; for $l = 0$ we take $\lambda_0 := v$ and fix $\{E_{0,j} : j \ge 0\}$ to be any enumeration of $\{E \in \FE(\L) : r(E) = v\}$. These satisfy~(i) by definition of $H_I$, and trivially satisfy (ii)~and~(iii).

Now as an inductive hypothesis, suppose that $l \ge 0$ and that $\lambda_1, \dots, \lambda_l$ and $\{E_{1,j} : j \ge 1\}, \dots, \{E_{l,j} : j \ge 1\}$ have been chosen and satisfy (i)--(iii). Just as in the proof of \cite[Lemma~4.7(1)]{Si1}, we have that $l \ge i_l$ so that $E_{i_l, j_l}$ has already been defined. If $\lambda_l(d(\lambda_{i_{l+1}}, d(\lambda_l))) \in E_{i_{l+1}, j_{l+1}}$ already, then $l > 0$ because $E \in \FE(\L)$ implies $E \cap \L^0 = \emptyset$, so $\lambda_{l+1} := \lambda_l$ and $E_{l+1, j} := E_{l,j}$ for all $j$ satisfy (i)--(iii) by the inductive hypothesis. On the other hand, if $\lambda_l(d(\lambda_{i_{l+1}}, d(\lambda_l))) \not\in E_{i_{l+1}, j_{l+1}}$, then $E := \Ext\big(\lambda_l(d(\lambda_{i_{l+1}}, d(\lambda_l))); E_{i_{l+1}, j_{l+1}}\big) \in \FE(\L)$ by \cite[Lemma~C.5]{RSY2}. By~\eqref{eq:saturation consequence}, there exists $\nu_{l+1} \in E$ such that $s(\nu) \not\in H_i$. But now $\lambda_{l+1} := \lambda_l\nu_{l+1}$ satisfies~(i) by choice of $\nu_{l+1}$, and taking $\{E_{l+1, j} : j \ge 1\}$ to be any enumeration of $\{E \in \FE(\L) : r(E) = s(\nu_{l+1})\}$ we have (ii)~and~(iii) satisfied just as in the proof of \cite[Lemma~4.7(1)]{Si1}.

The remainder of the proof of \cite[Lemma~4.7(1)]{Si1} shows that $x(0, d(\lambda_l)) := \lambda_l$ for all $l$ defines an element of $v\partial\L$, and since $H_I$ is hereditary, condition~(i) shows that $x(n) \not\in H_I$ for all $n \le d(x)$.
\hfil\penalty100\hbox{}\nobreak\hfill\hbox{\qed\ Claim~\ref{clm:boundary path}}
\renewcommand\qed{}\end{proof}

Now fix $w \in \L^0$. Let $x \in v\partial\L$ with $x(n) \not\in H_I$ for all $n$ as in Claim~\ref{clm:boundary path}. Since $\L$ is cofinal, there exists $n \le d(x)$ such that $w \L x(n) \not= \emptyset$. Since $x(n) \not\in H_I$ by construction of $x$, and since $H_I$ is hereditary, it follows that $w \not\in H_I$. Consequently $H_I = \emptyset$ as required.

Now suppose that $C^*(\L)$ is simple. Let $x \in \partial\L $, and let 
\[
H_x := \{w \in \L^0 : w \L x(n) = \emptyset\text{ for all $n$}\}.
\]
It is clear that $H_x$ is hereditary.  We claim that $H_x$ is saturated: suppose that $E \in v\FE(\L)$ with $s(E) \in H_x$, and suppose for contradiction that $\lambda \in v \L x(n)$.  If $\lambda = \mu\mu'$ for $\mu \in E$, then $\mu' \in s(\mu)\L x(n)$, contradicting $s(\mu) \in H_x$.  On the other hand, if $\lambda \not\in E\L$, then $\Ext(\lambda;E)$ is exhaustive by \cite[Lemma~2.3]{Si1}.  Since $x \in \partial(\L;\Ee)$, it follows that $x(n, n + d(\alpha)) = \alpha$ for some $\alpha \in \Ext(\lambda;E)$; say $(\alpha,\beta) \in \Lmin(\lambda,\mu)$ where $\mu \in E$.  Then $\beta \in s(\mu)\L x(n+d(\alpha))$, again contradicting $s(\mu) \in H_x$.  This proves our claim.

Now $H_x \not = \L^0$ because, in particular, $r(x) \not \in H_x$.  It follows that if $H_x$ is nonempty then it corresponds to a nontrivial ideal $I_{H_x}$ which is impossible since $C^*(\L)$ is simple by assumption.  Hence $\L$ is cofinal as required.
\end{proof}

\begin{dfn} \label{dfn:loop with entrance}
Let $(\L,d)$ be a finitely aligned $k$-graph. We say that a path $\mu \in \L$ is a \emph{loop with an entrance} if $s(\mu) = r(\mu)$ and there exists $\alpha \in s(\mu)\L$ such that $d(\mu) \ge d(\alpha)$ and $\mu(0, d(\alpha)) \not= \alpha$. We say that a vertex $v \in \L^0$ can be \emph{reached from a loop with an entrance} if there exists a loop with an entrance $\mu \in \L$ such that $v \L s(\mu) \not= \emptyset$. 
\end{dfn}

The following proposition rectifies a slight error in \cite[Proposition~4.9]{KP}, specifically in the argument that $\mathcal{G}_\Lambda$ is locally contracting. Our condition that every vertex can be reached from a loop with an entrance is slightly than that in \cite{KP} that every vertex can be reached from a nontrivial loop, and this stronger condition is needed to make both our argument and that of \cite{KP} run.

\begin{prop}\label{prp:purely infinite} 
Let $(\L,d)$ be a finitely aligned $k$-graph, and suppose that $\L$ satisfies condition~{\rm(C)}.  Suppose also that every $v \in \L^0$ can be reached from a loop with an entrance. Then every nontrivial hereditary subalgebra of $C^*(\L)$ contains an infinite projection. In particular, if $\L$ is also cofinal, then $C^*(\L)$ is purely infinite.
\end{prop}

The proof of Proposition~\ref{prp:purely infinite} is based heavily on the proof of \cite[Proposition~5.3]{BPRS}. First we need to recall some definitions and establish some technical results and notation. Definitions \ref{dfn:Ts and nus}~and~\ref{dfn:Psub n,v} and the proof of Lemma~\ref{lem:norm equality} are based almost entirely on the definitions and techniques used in \cite{RSY2} from \cite[Notation~3.12]{RSY2} to the proof of \cite[Proposition~3.13]{RSY2}. We present them seperately here because the conclusion of Lemma~\ref{lem:norm equality} is not stated explicitly in \cite{RSY2}.

\begin{dfn}\label{dfn:Ts and nus}
Let $(\L,d)$ be a finitely aligned $k$-graph, and let $E \subset \L$ be finite. As in \cite[Notation~3.12]{RSY2}, for all $n$ and $v$ such that $(\Pi E) v \cap \L^n$ is nonempty, we write $T^{\Pi E}(n,v)$ for the set $\{\nu \in v\L \setminus \{v\} : \lambda\nu \in \Pi E\text{ for some }\lambda \in (\Pi E) v \cap \L^n\}$. By the properties of $\Pi E$, the set $T(\lambda) := \{\nu \in s(\lambda)\L\setminus\{s(\lambda)\} : \lambda\nu \in \Pi E\}$ is equal to $T^{\Pi E}(n,v)$ for all $\lambda \in (\Pi E)v \cap \L^n$ \cite[Remark~3.4]{RSY2}. If, in addition to $(\Pi E)v \cap \L^n \not= \emptyset$, we have $T^{\Pi E}(n,v) \not\in \FE(\L)$, we fix, once and for all, an element $\xi^{\Pi E}(n,v)$ of $v \L$ such that $\Ext(\xi^{\Pi E}(n,v); T^{\Pi E}(n,v)) = \emptyset$, and for $\lambda \in (\Pi E)v \cap \L^n$, we define $\xi_\lambda := \xi^{\Pi E}(n,v)$.

Notice that if $\lambda,\mu \in \Pi E$ satisfy $s(\lambda) = s(\mu)$ and $d(\lambda) = d(\mu)$, then we also have $T(\lambda) = T(\mu)$ and $\xi_\lambda = \xi_\mu$.
\end{dfn}

\begin{dfn}\label{dfn:Psub n,v}
Let $(\L,d)$ be a finitely aligned $k$-graph, let $E \subset \L$ be finite, and let $\{t_\lambda : \lambda \in \L\}$ be a Cuntz-Krieger $\L$-family. For each $n, v$ such that $(\Pi E)v \cap \L^n$ is nonempty and $T^{\Pi E}(n,v)$ is not exhaustive, we define
\[
P_{n,v} := \sum_{\lambda\in (\Pi E)v \cap \L^n} s_{\lambda{\xi_\lambda}} s^*_{\lambda{\xi_\lambda}} \in C^*(\L).
\]
\end{dfn}

\begin{notation}
Let $(\L,d)$ be a finitely aligned $k$-graph. We write $\Phi$ for the linear map from $C^*(\L)$ to $C^*(\L)^\gamma$ determined by $\Phi(a) := \int_\TT \gamma_z(a)\,dz$. We have that $\Phi$ is positive and is faithful on positive elements.
\end{notation}

\begin{lemma}\label{lem:norm equality} 
Let $(\L,d)$ be a finitely aligned $k$-graph, let $E \subset \L$ be finite, and let $a = \sum_{\lambda,\mu \in \Pi E} a_{\lambda,\mu} s_\lambda s^*_\mu$ with $a \not= 0$. For $n \in \NN^k$ and $v \in \L^0$ such that $(\Pi E)v \cap \L^n$ is nonempty and $T^{\Pi E}(n,v)$ is not exhaustive, let
\[
\Ff_{\Pi E}(n,v) := \clsp\{s_{\lambda{\xi_\lambda}} s^*_{\mu{\xi_\lambda}} : \lambda,\mu \in (\Pi E) v \cap \L^n\}.
\]
Then for all $n,v$ such that $(\Pi E) v \cap \L^n$ is nonempty and $T^{\Pi E}(n,v)$ is not exhaustive, we have that $P_{n,v} \Phi(a) \in \Ff_{\Pi E}(n,v)$.  Furthermore, there exist $n_0, v_0$ such that $(\Pi E) v_0 \cap \L^{n_0}$ is nonempty and $T^{\Pi E}(n_0, v_0)$ is not exhaustive, and such that $\|P_{n_0, v_0} \Phi(a)\| = \|\Phi(a)\|$.
\end{lemma}
\begin{proof}
By \cite[Lemma~3.15]{RSY2}, we have that each $s_{\lambda\xi_\lambda} s^*_{\lambda\xi_\lambda} \le Q(s)^{\Pi E}_\lambda$ where $Q(s)^{\Pi E}_\lambda$ is defined by \eqref{eq:Qdef}. Since the $Q(s)^{\Pi E}_\lambda$ are mutually orthogonal projections, it follows that $s_{\lambda\xi_\lambda} s^*_{\lambda\xi_\lambda} Q(s)^{\Pi E}_\mu = \delta_{\lambda,\mu} s_{\lambda\xi_\lambda} s^*_{\lambda\xi_\lambda}$. Hence, for $(\lambda,\mu) \in \Pi E \times_{d,s} \Pi E$, we have
\begin{equation}\label{eq:left mpctn}
P_{n,v} \Theta(s)^{\Pi E}_{\lambda,\mu} 
= P_{n,v} Q(s)^{\Pi E}_\lambda s_\lambda s^*_\mu 
= s_{\lambda\xi_\lambda} s^*_{\lambda\xi_\lambda} s_\lambda s^*_\mu
= s_{\lambda\xi_\lambda} s^*_{\mu\xi_\lambda},
\end{equation}
and hence $P_{n,v} \Phi(a) \in \Ff_{\Pi E}(n,v)$. Moreover, taking adjoints in~\eqref{eq:left mpctn}, shows that each $P_{n,v}$ commutes with each $\Theta(s)^{\Pi E}_{\lambda,\mu}$.

By definition of the $\Theta(s)^{\Pi E}_{\lambda,\mu}$, and by \cite[Corollary~4.10]{Si1}, we have that $\Theta(s)^{\Pi E}_{\lambda,\mu}$ is nonzero if and only if $T(\lambda)$ is not exhaustive. Moreover, since the $Q(s)^{\Pi E}_\lambda$ are mutually orthogonal and dominate the $s_{\lambda\xi_\lambda} s^*_{\lambda\xi_\lambda}$, we have that the latter are also mutually orthogonal. It follows from this and from~\eqref{eq:left mpctn} that 
\[
b \mapsto \sum_{\substack{(\Pi E)v \cap \L^n \not= \emptyset \\ T^{\Pi E}(n,v) \not\in \FE(\L)}} P_{n,v} b
\]
is an injective homomorphism of $\clsp\{\Theta(s)^{\Pi E}_{\lambda,\mu} : \lambda,\mu \in \Pi E \times_{d,s} \Pi E\}$. Since injective $C^*$-homomorphisms are isometric, it follows that $\big\|\sum P_{n,v} \Phi(a)\big\| = \|\Phi(a)\|$.

Since the $P_{n,v}$ are mutually orthogonal and commute with $\Phi(a)$, there therefore exists a vertex $v_0$ and a degree $n_0$ such that $\|\Phi(a)\| = \|P_{n_0, v_0} \Phi(a)\|$. Clearly for this $n_0,v_0$ we must have $(\Pi E)v_0 \cap \L^{n_0}$ nonempty and $T(\lambda)$ non-exhaustive for $\lambda \in (\Pi E)v_0 \cap \L^{n_0}$, for otherwise we have $P_{n_0, v_0} = 0$ contradicting $a \not= 0$.
\end{proof}

\begin{lemma}\label{lem:infinite vert projs} 
Let $(\L,d)$ be a finitely aligned $k$-graph, and suppose that every $v \in \L^0$ can be reached from a loop with an entrance. Then for each $v \in \L^0$, the projection $s_v$ is infinite, and hence for each $\lambda \in \L$, the range projection $s_\lambda s^*_\lambda$ is also infinite.
\end{lemma}
\begin{proof} 
Fix $v \in \L^0$, and let $\mu$ be a loop with an entrance such that $v\L s(\mu)$ is nonempty. Fix $\lambda \in v\L s(\mu)$, and fix $\alpha \in s(\mu)\L$ such that $d(\alpha) \le d(\mu)$ and $\mu(0, d(\alpha)) \not= \alpha$. We have $s_v \ge s_\lambda s^*_\lambda \sim s^*_\lambda s_\lambda = s_{s(\mu)}$, so it suffices to show that $s_{s(\mu)}$ is infinite. But (TCK3) ensures that $s_\mu s^*_\mu s_\alpha s^*_\alpha = 0$, and it follows that $s_{s(\mu)} = s^*_\mu s_\mu \sim s_\mu s^*_\mu \le s_{s(\mu)} - s_\alpha s^*_\alpha < s_{s(\mu)}$. 

For the last statement, notice that $s_{s(\lambda)}$ is infinite by the previous paragraph, and $s_\lambda s^*_\lambda \sim s^*_\lambda s_\lambda = s_{s(\lambda)}$.
\end{proof}

\begin{lemma}[{\cite[Lemma~5.4]{BPRS}}] \label{lem:BPRS lem} 
Let $E \subset \L^n$, let $w \in s(E)$, and let $t$ be a positive element of $\Ff_E(w) := \lsp\{s_\lambda s^*_\mu : \lambda,\mu \in Ew\}$.  Then there is a projection $r$ in $C^*(t) \subset \Ff_E(w)$ such that $r t r = \|t\| r$.
\end{lemma}
\begin{proof} 
The proof is formally identical to that of \cite[Lemma~5.4]{BPRS}
\end{proof}

\begin{proof}[Proof of Proposition~\ref{prp:purely infinite}] 
Our proof follows that of \cite[Proposition~5.3]{BPRS} very closely.

Fix a nontrivial hereditary subalgebra $A$ of $C^*(\L)$, and a positive element $a \in A$ such that $\Phi(a) \in C^*(\L)^\gamma$ satisfies $\|\Phi(a)\| = 1$.  Let $b = \sum_{\lambda,\mu \in E} b_{\lambda,\mu} s_\lambda s^*_\mu$ be a finite linear combination such that $b > 0$ and $\|a - b\| \le \frac{1}{4}$; this is always possible because $\lsp\{s_\lambda s^*_\mu : \lambda,\mu \in \L\}$ is a dense $^*$-subalgebra of $C^*(\L)$.  Let $b_0 := \Phi(b)$.  Since $\Phi$ is norm-decreasing and linear, we have
\[
1 - \|b_0\| = \big|\|\Phi(a)\| - \|\Phi(b)\|\big| \le \|\Phi(a - b)\| \le \|a-b\| \le \frac{1}{4},
\]
and hence $\|b_0\| \ge \frac{3}{4}$.  Furthermore, $b_0 \ge 0$ because $\Phi$ is positive.  Applying Lemma~\ref{lem:norm equality}, we obtain a projection $P_{n_0, v_0}$ such that $b_1 := P_{n_0, v_0} b_0$ satisfies $b_1 \in \Ff_{\Pi E}(n_0,v_0)$ and $\|b_1\| = \|b_0\|$, where $(\Pi E)v_0 \cap \L^{n_0}$ is nonempty and $T^{\Pi E}(n_0, v_0)$ is not exhaustive.  Notice that $b_1 \ge 0$.  By Lemma~\ref{lem:BPRS lem} there exists a projection $r \in C^*(b_1)$ with $r b_1 r = \|b_1\| r$; note that $r$ is clearly nonzero.  Let $v_1 := s(\xi^{\Pi E}(n_0, v_0))$, and let $S := \{\lambda{\xi_\lambda} : \lambda \in (\Pi E) v_0 \cap \L^{n_0}\}$.

Since $b_1 \in \lsp\{s_{\lambda} s^*_{\mu} : \lambda,\mu \in S\}$, which is a matrix algebra indexed by $S$, we can express $r$ as a finite sum $r = \sum_{\lambda,\mu \in S} r_{\lambda,\mu} s_\lambda s^*_\mu$, and the $S \times S$ matrix $(r_{\lambda,\mu})$ is a projection.

Since $(\L,d)$ satisfies condition~(C), there exists $x \in v_1 \partial\L$ such that for $\lambda,\mu \in \L r(x)$ with $\lambda \not= \mu$, we have $\MCE(\lambda x, \mu x) = \emptyset$. By \cite[Lemma~6.4]{Si1}, for distinct $\lambda,\mu \in S$, there exists $n^x_{\lambda,\mu}$ such that $\Lmin(\lambda x(0, n^x_{\lambda,\mu}), \mu x(0, n^x_{\lambda,\mu})) = \emptyset$. Let
\[
\textstyle 
M := \bigvee\{n^x_{\lambda,\mu} : \lambda,\mu \in S, \lambda\not=\mu\},
\]
and let $x_M := x(0,M)$.  Let $q := \sum_{\lambda,\mu \in S} r_{\lambda,\mu} s_{\lambda x_M} s^*_{\mu x_M}$.  Since the matrix $(r_{\lambda,\mu})$ is a nonzero projection in $M_S(\CC)$, we know that $q$ is a nonzero projection in $\Ff_{N_E + d(x_M)}$, and since $s_{x_M} s^*_{x_M}$ is a subprojection of $s_{v_1}$, we have $q \le r$. Using the defining property of $x_M$ as in the proof of \cite[Lemma~6.7]{Si1}, we have that $q P_{n_0, v_0} b q = q P_{n_0, v_0} b_0 q = q b_1 q$. Now $q \le P_{n_0, v_0}$ by definition so our choice of $r$ gives
\[
q b q = q b_1 q = q r b_1 r q = \|b_1\|rq = \|b_0\| q \ge \frac{3}{4} q.
\]
Since $\|a - b\| \le \frac{1}{4}$, we have $qaq \ge qbq - \frac{1}{4}q \ge \frac{3}{4} q - \frac{1}{4}q = \frac{1}{2} q$, and it follows that $q a q$ is invertible in $q C^*(\L) q$.  Write $c$ for the inverse of $q a q$ in $q C^*(\L) q$, and let
\[
t := c^{1/2} q a^{1/2}.  
\]
Then $t^* t = a^{1/2} q c q a^{1/2} \le \|c\| a$, so $t^* t \in A$ because $A$ is hereditary.

We now need only show that $t^* t$ is an infinite projection. But 
\[
t^* t \sim t t^* = c^{1/2} q a q c^{1/2} = 1_{q C^*(\L) q} = q,
\]
so it suffices to show that $q$ is infinite. By choice of $n_0, v_0$, there exists $\sigma \in S$.  By Lemma~\ref{lem:infinite vert projs}, $s_{\sigma x_M} s^*_{\sigma x_M}$ is infinite. But $s_{\sigma x_M} s^*_{\sigma x_M}$ is a minimal projection in the finite-dimensional $C^*$-algebra $\lsp\{s_{\sigma x_M} s^*_{\tau x_M} : \sigma,\tau \in S\}$, which contains $q$. Since $q \not= 0$, $s_{\sigma x_M} s^*_{\sigma x_M}$ is equivalent to a subprojection of $q$, so $q$ is infinite.
\end{proof}

\begin{cor}
Let $(\L,d)$ be a finitely aligned $k$-graph. Suppose that $\L$ satisfies condition~{\rm(C)} and is cofinal, and that every $v \in \L^0$ can be reached from a loop with an entrance. Then $C^*(\L)$ is determined up to isomorphism by its $K$-theory.
\end{cor}
\begin{proof}
We have that $C^*(\L)$ is nuclear and satisfies UCT by Proposition~\ref{prp:rel algs nuclear}, is simple by Proposition~\ref{prp:simple graph alg}, and is purely infinite by Proposition~\ref{prp:purely infinite}. The result then follows from the Kirchberg-Phillips classification theorem \cite[Theorem~4.2.4]{P}.
\end{proof}


\begin{thebibliography}{00}
\bibitem{BPRS} T. Bates, D. Pask, I. Raeburn, and W. Szyma\'nski, \emph{The $C^*$-algebras of row--finite graphs}, New York J. Math.  {\bf 6} (2000), 307--324.

\bibitem{BHRS} T. Bates, J. Hong, I. Raeburn, and W. Szyma\'nski, \emph{The ideal structure of the $C^*$-algebras of infinite graphs}, Illinois J. Math. {\bf 46} (2002), 1159--1176.

\bibitem{HS} J. H. Hong and W. Szyma\'nski, \emph{The primitive ideal space of the $C^*$-algebras of infinite graphs}, J. Math. Soc. Japan {\bf 56} (2004), 45--64.

\bibitem{KP} A. Kumjian and D. Pask, \emph{Higher rank graph $C^*$-algebras}, New York J. Math. {\bf 6} (2000), 1--20.

\bibitem{KPR} A. Kumjian, D. Pask, and I. Raeburn, \emph{Cuntz-Krieger algebras of directed graphs,} Pacific J. Math {\bf 184} (1998), 161--174.

\bibitem{PQR} D. Pask, J.C. Quigg, and I. Raeburn, \emph{Coverings of $k$-graphs}, preprint, 2004 [arXiv:math.OA/\allowbreak0401017].

\bibitem{P} N.C. Phillips, \emph{A classification theorem for nuclear purely infinite simple $C^*$-algebras}, Documenta Math. {\bf 5} (2000), 49--114.

\bibitem{Qui} J.C. Quigg, \emph{Discrete coactions and $C^*$-algebraic bundles}, J. Austral.  Math.  Soc (Series A) {\bf 60} (1996), 204--221.

\bibitem{RS1} I. Raeburn and A. Sims, \emph{Product systems of graphs and the Toeplitz algebras of higher-rank graphs}, J. Operator Th., to appear [arXiv:math.OA/0305371].

\bibitem{RSY1} I. Raeburn, A. Sims and T. Yeend \emph{Higher-rank graphs and their $C^*$-algebras}, Proc. Edinb. Math. Soc. {\bf 46} (2003), 99--115.

\bibitem{RSY2} I. Raeburn, A. Sims and T. Yeend \emph{The $C^*$-algebras of finitely aligned higher-rank graphs}, J. Funct. Anal. {\bf213} (2004), 206--240.

\bibitem{RSc} J. Rosenberg and C. Schochet, \emph{The K\"unneth theorem and the universal coefficient theorem for Kasparov's generalised $K$-functor}, Duke Math. J. {\bf 55} (1987), 431--474.

\bibitem{Si1} A. Sims, \emph{Relative Cuntz-Krieger algebras of finitely aligned higher-rank graphs}, Indiana U. Math. J., to appear [arXiv:math.OA/0312152].

\bibitem{Sz2} W. Szyma\'nski, \emph{Simplicity of Cuntz-Krieger algebras of infinite matrices}, Pac.  J. Math.  {\bf 199} (2001), 249--256.

\end{thebibliography}
\end{document}